\documentclass[letter,fleqn,12pt]{article}
\usepackage{graphicx}
\usepackage{algorithmic,algorithm}

\usepackage{bbold, epsfig}
\usepackage{latexsym, amsmath, amssymb, euscript}

\usepackage{fontenc, indentfirst}
\usepackage{fancyhdr}
\usepackage{multicol}
\usepackage{subfig}
\usepackage{placeins}
\usepackage[usenames,dvipsnames]{xcolor}
\usepackage[english]{babel}
\usepackage{textcomp}
\usepackage{bbm}

\usepackage{amsmath}
\usepackage{amssymb}
\usepackage{theorem}
\usepackage{ifthen}

\newcommand{\revt}[1]{{\color{black}#1}}
\newcommand{\revtt}[1]{{\color{black}#1}}

\newtheorem{theorem}{Theorem}[section] 
 
\newtheorem{proposition}[theorem]{Proposition} 
\newtheorem{definition}[theorem]{Definition} 
\newtheorem{remark}[theorem]{Remark} 
\newtheorem{corollary}[theorem]{Corollary}

\newcommand{\mR}{\mathbb R}
\newcommand{\mN}{\mathbb N}

\newcommand{\cO}{\mathcal O}

\def\hu{{\widehat{u}}}
\def\hlambda{{\widehat{\lambda}}}

\def\heta{{\widehat{\eta}}}

\def\Dpartial#1#2{ {\frac{\partial #1}{\partial #2} }}
\def\Dpartialmix#1#2#3{ {\partial^2 #1 \over \partial #2\,\partial #3} }

\DeclareMathOperator{\rank}{rank}
\DeclareMathOperator{\Span}{span}
\DeclareMathOperator{\argmin}{argmin}

\begin{document}
\thispagestyle{empty}




%
%


\setcounter{page}{1}
\begin{center}
  \textbf{\Large {A gradient method in a Hilbert space with an
      optimized inner product:} achieving a Newton-like convergence }
\end{center}

\begin{center}
Arian {\sc Novruzi}\footnote{\revt{Department of Mathematics and Statistics, University of Ottawa, Ottawa, ON, K1N 6N5, Canada, Email: {\tt novruzi@uottawa.ca}; the corresponding author}} and
Bartosz {\sc Protas}\footnote{\revt{Department of Mathematics and Statistics, McMaster University, Hamilton, ON, L8S 4K1, Canada, Email: {\tt bprotas@mcmaster.ca}}}
\end{center}

\begin{abstract}
  In this paper we introduce a new gradient method which attains
  quadratic convergence in a certain sense.  Applicable to
  infinite-dimensional {unconstrained minimization problems posed in a
    Hilbert space $H$}, the approach consists in finding the energy
  gradient $g(\lambda)$ {defined with respect to an optimal} inner
  product {selected from an infinite family of equivalent} inner
  products $(\cdot,\cdot)_\lambda$ in the space $H$. {The inner
    products are} parameterized by a space-dependent weight \revt{function}
  $\lambda$.
  At each iteration of the method, \revt{where an approximation to the
    minimizer is given by an element $u\in H$,} an optimal weight
  $\hlambda$ is found as a solution of a nonlinear minimization
  problem in the space of weights $\Lambda$.  It turns out that the
  projection of $\kappa g(\hlambda)$, where $0<\kappa \ll 1$ is a
  fixed step size, onto a certain finite-dimensional subspace
  generated by the method is consistent with Newton's step $h$, in the
  sense that $P_u(\kappa g(\hlambda))=P_u(h)$, where $P_u$ is \revt{an
    operator describing the projection onto the subspace}. As
  demonstrated by rigorous analysis, this property ensures that thus
  constructed gradient method attains quadratic convergence for error
  components contained in these subspaces, in addition to the linear
  convergence typical of the standard gradient method. We propose a
  numerical implementation of {this} new approach and analyze its
  complexity.  Computational results obtained based on a simple model
  problem confirm the theoretically established convergence
  properties, demonstrating that the proposed approach performs much
  better than the standard steepest-descent method based on Sobolev
  gradients. The presented results offer an explanation of a number of
  earlier empirical observations concerning the convergence of
  Sobolev-gradient methods.
\end{abstract}

\tableofcontents

\section{Introduction}

In this investigation we consider 
solution of general {unconstrained} optimization problems using the
steepest-descent method and focus on {modifying the definition of the
  gradient} such that in certain circumstances {this approach} will
achieve a quadratic convergence, characteristic of Newton's method.
This is accomplished by judiciously exploiting the freedom inherent in
the choice of different equivalent norms defining the gradient through
the Riesz representation theorem. This freedom can be used to adjust
the definition of the inner product, such that the resulting gradient
will, in a suitable sense, best resemble the corresponding Newton
step.  While such ideas can be pursued in both finite-dimensional and
infinite-dimensional settings, the formulation is arguably more
interesting mathematically and more useful in practice in the latter
case. We will thus consider unconstrained optimization problems of the
general form
\begin{equation}
\label{e:min{e(u)}}
e(\hu) = \min\{e(u),\; u \in H \},
\end{equation}
where $e \; : \; H \mapsto \mR$ is the objective functional and $H$ is
a suitable function space with Hilbert structure. 
Applications which have this form include, for example,
minimization of various energy functionals in physics and optimization
formulations of inverse problems, where evaluation of the objective
functional $e(u)$ may involve solution of a complicated (time-)dependent
partial differential equation (PDE). In such applications $H$ is typically
taken as a Sobolev space $H^p(\Omega)$, $p \in \mathbb N$, where $\Omega
\subset \mR^d$ is the spatial domain assumed to be sufficiently smooth (Lipschitz)
and $d \ge 1$ is the dimension \cite{af05}.  Therefore, to fix
attention, here we will assume $H := H^1_0(\Omega)$ with the inner
product and norm in $H$ defined as
\begin{equation}
\label{e:(u,v)_Ht}
(u,v)_H=\int_{\Omega}((\nabla u\cdot\nabla v) + uv)\, dx,
\qquad
\|u\|_H=(u,v)_H^{1/2}.
\end{equation}

\begin{algorithm}
\begin{algorithmic}[1]
\STATE set $u = u_0$
\REPEAT 
\STATE  compute $e'(u;\cdot)$ 
\STATE  determine $g\in H$ as the solution of $(g,v)_H=e'(u;v)$, $\forall v\in H$
\STATE  set $\widetilde{u} = u - \kappa g$
\STATE  set $u = \widetilde{u}$
\UNTIL{ \  $|e(u)|<$ \texttt{err}}
\end{algorithmic}
\caption{
 (Sobolev) Gradient method  \newline
     \textbf{Input:} \newline
 \hspace*{0.5cm} $u_0 \in H$ --- initial guess \newline
 \hspace*{0.5cm} $\kappa > 0$ --- step size (sufficiently small) \newline
 \hspace*{0.5cm} $\texttt{err} > 0$ --- tolerance \newline
 \textbf{Output:} \newline
 \hspace*{0.5cm} $\widetilde{u}$ --- approximation to the solution $\hu$ of problem \eqref{e:min{e(u)}}}
\label{alg:gradient}
\end{algorithm}

\begin{algorithm}
\begin{algorithmic}[1]
\STATE set $u = u_0$
\REPEAT 
\STATE  {compute $e'(u;\cdot)$}
\STATE  {compute $e''(u;\cdot,\cdot)$}
\STATE  determine $h\in H$ as the solution of $e''(u;h,v)=e'(u;v)$, $\forall v\in H$
\STATE  set $\widetilde{u} = u - h$
\STATE  set $u = \widetilde{u}$
\UNTIL{ \  $|e(u)|<$ \texttt{err}}
\end{algorithmic}
\caption{
 Newton's method  \newline
     \textbf{Input:} \newline
 \hspace*{0.5cm} $u_0 \in H$ --- initial guess \newline
 \hspace*{0.5cm} $\texttt{err} > 0$ --- tolerance \newline
 \textbf{Output:} \newline
 \hspace*{0.5cm} $\widetilde{u}$ --- approximation to the solution $\hu$ of problem \eqref{e:min{e(u)}}}
\label{alg:newton}
\end{algorithm}

The two most elementary approaches to solve problem
\eqref{e:min{e(u)}} are the gradient and Newton's method which, for
the sake of completeness, are defined in Algorithms \ref{alg:gradient}
and \ref{alg:newton}, respectively. The former approach is sometimes
also referred to as the ``Sobolev gradient'' method \cite{n10}. While
gradient approaches often involve an adaptive step size selection
\cite{nw00}, for simplicity of analysis in Algorithm
\ref{alg:gradient} we consider a fixed step size $\kappa =
\text{Const}$. Likewise, in order to keep the analysis tractable, we
do not consider common modifications of the gradient approach such as
the conjugate-gradient method.  As regards convergence of the gradient
and Newton's method, we have the following classical results, see for
example \cite{ciarlet-1}.
\begin{theorem}
\label{th:gradient}
Let $\hu\in H$ be {a} solution of (\ref{e:min{e(u)}}). Assume
that $e$ is twice differentiable near $\hu$ and there exist
$\delta_0>0$ and $\alpha_0>0$ such that $e''(u;v,v)\geq
\alpha_0\|v\|^2_H$ for all $\|u-\hu\|_H\leq \delta_0$ and $v\in H$.  Then,
for all $u_0\in B_H(\hu,\delta_0)$ the gradient {method given by}
Algorithm \ref{alg:gradient} converges linearly to $\hu$ in $H$, i.e.,
\begin{eqnarray*}
\|\widetilde{u}-\hu\|_H
&\leq&
\epsilon\|u-\hu\|_H,
\end{eqnarray*}
with $\epsilon\in(0,1)$  depending only on $e$ and $\kappa$.
\end{theorem}

\begin{theorem}\label{th:newton}
  Let $\hu\in H$ be {a} solution of (\ref{e:min{e(u)}}). Assume
  that $e$ is three times differentiable near $\hu$ and that the map
  $v\in H\mapsto e''(\hu;v,\cdot)\in H'$ is invertible.  Then, there
  exists $\delta_0>0$ such that for all $u_0\in B_H(\hu,\delta_0)$
  Newton's {method given by} Algorithm \ref{alg:newton} converges
  quadratically to $\hu$ in $H$, i.e.,
\begin{eqnarray*}
\|\widetilde{u}-\hu\|_H
&\leq& C \|u-\hu\|^2_H	,
\end{eqnarray*}
with $C>0$ depending only on $e$.
\end{theorem}

We emphasize that in Algorithm \ref{alg:gradient} the gradient $g$
must be computed with respect to the topology of the space $H$ in
which the solution {to \eqref{e:min{e(u)}}} is sought \cite{n10},
an aspect of the problem often neglected in numerical investigations.
A metric equivalent to $\| \cdot \|_H$ (in the precise sense of norm
equivalence) can be obtained by redefining the inner product in
\eqref{e:(u,v)_Ht} more generally as follows
\begin{equation}
\label{e:(u,v)_lam0}
(u,v)_{\lambda_0} 
= 
\int_{\Omega}(\lambda_0 (\nabla u\cdot\nabla v) + uv)\, dx,
\end{equation}
where $\lambda_0\in(0,\infty)$ is a fixed constant. While as compared
to \eqref{e:(u,v)_Ht} definition \eqref{e:(u,v)_lam0} does not change
the analytic structure of the optimization problem
\eqref{e:min{e(u)}}, there is abundant computational evidence obtained
in the solution of complicated optimization problems
\cite{pbh02,p07b,ap15a} that convergence of gradient Algorithm
\ref{alg:gradient} may be significantly accelerated by replacing the
inner product from \eqref{e:(u,v)_Ht} with the one introduced in
\eqref{e:(u,v)_lam0} for some judicious {choices} of the
parameter $\lambda_0$. Likewise, a similar acceleration was also
observed when the inner product in \eqref{e:(u,v)_Ht} was replaced
with another equivalent definition motivated by the structure of the
minimization problem and different from \eqref{e:(u,v)_lam0},
cf.~\cite{rssl09,ms10,kd12}. In the absence of an understanding of the
mechanism responsible for this acceleration, the parameter
$\lambda_0$, or other quantities parameterizing the equivalent inner
product, were chosen empirically by trial and error, which is
unsatisfactory. 

In the present investigation we will consider a more general form of
the inner product \eqref{e:(u,v)_lam0} in which the constant
$\lambda_0$ is replaced with a space-dependent weight $\lambda =
\lambda(x)$. Our goal is to develop a rational and systematic approach
allowing one to accelerate the convergence of gradient iterations in
Algorithm \ref{alg:gradient} in comparison to the standard case by
adaptively adjusting the weight $\lambda(x)$. This will result in a
reduction of {the total number of iterations needed to solve
  problem \eqref{e:min{e(u)}} to a given accuracy,} but each iteration
will be more costly.

Modifications of the inner product with respect to which the gradient
is defined may also be interpreted as gradient preconditioning and
this perspective is pursued in the monograph \cite{fk02} focused on
related problems arising in the solution of nonlinear elliptic
equations. The relationship between the gradient and Newton's methods
was explored in \cite{kn07} where a variable inner product was
considered.  In contrast to \revt{the present} approach \revt{in which
  the} inner-product \revt{weights are} sought by matching the
projections of \revt{the} gradient and Newton's steps \revt{onto} a
certain subspace, \revt{in \cite{kn07} optimal inner products were
  found by maximizing the descent achieved at a given iteration with
  respect to the structure of the corresponding preconditioning
  operator.}

The structure of the paper is as follows.  In the next section we
define the new approach in a general form, whereas in Section
\ref{sec:errors} we prove its convergence properties. Then, in Section
\ref{sec:optimal} we describe the numerical approach implementing the
general method introduced in Section \ref{sec:newgrad} {in two
  variants} and analyze its computational complexity.  Our model
problem and computational results are presented in Section
\ref{sec:results}, whereas discussion and conclusions are deferred to
Section \ref{sec:final}.

\section{A new gradient method based on an optimal inner product}
\label{sec:newgrad}

In this section we introduce our modified version of Algorithm
\ref{alg:gradient} for the solution of the minimization problem
\eqref{e:min{e(u)}}. We begin by making the following assumptions on
the energy $e(u)$:
\begin{eqnarray}
e&&\mbox{\it is $C^2$ in $H$},
\\
|e''(u;v,w)|
&\leq&
M\|u\|_H \|w\|_H,\quad u,v,w \in H,
\label{e:e->C0}\\
e''(u;v,v)
&\geq&
m\|v\|_H^2,
\label{e:e->convex}
\end{eqnarray}
with certain $M,m>0$.

Let us point out that at each step of both the gradient and Newton's
methods the descent direction is defined by the solution of the
equation
\begin{equation}\label{e:b=e'->1}
 b(z,v)=e'(u;v),\quad \forall v\in H,
\end{equation}
where $b(u,v)$, a symmetric bilinear continuous elliptic form, and
$z\in H$ are specific to each method, namely,
\begin{itemize}
\item $b(u,v)=(u,v)_H$ and $z=g$ in the case of {the} gradient method, and 
\item $b(u,v)=e''(u;u,v)$ and $z=h$ in the case of Newton's method.
\end{itemize}
Moreover, we note that the solution $z$ of \eqref{e:b=e'->1} is also
the solution of the minimization problem
\begin{equation}
\label{e:min{b-l}}
\min\left\{\frac{1}{2}b(v,v)-e'(u;v),\;\, v\in H\right\}. 
\end{equation}
We emphasize that, in fact, Newton's method may be also viewed as a
``gradient'' method with a particular choice of the inner product at
each iteration, namely, the one induced by $e''(u)$.  Therefore, the
idea for improving the classical gradient method is to make the
gradient step $g$ ``close'' to the Newton step $h$ by suitably
adapting the inner product in $H$.

We thus propose the following modification of the gradient method from
Algorithm \ref{alg:gradient}.  We want to consider the gradient $g$
{defined} with respect to {an} inner product in $H$
depending on a function parameter $\lambda$. Typically, $0<\lambda\in
C^0(\overline{\Omega})$, {however, to make our method more
  attractive from the computational point of view} we will consider
$\lambda$ with a finite range.  Namely, let 
$\{\Omega_{i},\, i=1,\ldots,N\}$, be a partition of $\Omega$ into open Lipschitz sets,
$\Lambda=\{\lambda:\Omega\mapsto\mR,\,\; \lambda(\Omega_i)=\lambda_i
\in\mR\} = \Span \{\ell^i,\, i=1,\dots,N\}\subset L^\infty(\Omega)$,
where $\ell^i\in\Lambda$, $\ell^i=\delta_{i,j}$ in $\Omega_j$,
$i,j=1,\ldots,N$ with $\delta_{i,j}$ the Kronecker symbol,
$\Lambda^+=\Lambda\cap\{0<\lambda<\infty\}$.  Sometimes {without
  the risk of confusion we will write
  $\lambda=[\lambda_1,\ldots,\lambda_N]\in\mathbb R^N$ for
  $\lambda\in\Lambda$}, meaning $\lambda(\Omega_i)=\lambda_i$, for all
$i=1,\ldots, N$.  Then, for $\lambda\in \Lambda^+$, we define the
following inner product and norm in $H$
\begin{equation}\label{e:(u,v)_lambda}
(v,w)_\lambda
=
\int_{\Omega} 
\lambda (\nabla v\cdot \nabla w) + v w \, dx,
\quad
\|v\|_\lambda=(v,v)_\lambda^{1/2},
\quad
\forall v,w\in H.
\end{equation}
Clearly, $(\cdot,\cdot)_\lambda$ and $(\cdot,\cdot)_H=(\cdot,\cdot)_1$ are
equivalent in $H$ and therefore we can use $(\cdot,\cdot)_\lambda$
instead of $(\cdot,\cdot)_H$ for the gradient method.

The idea is to use the inner product $(\cdot,\cdot)_\lambda$ in the
gradient method, with $\lambda$ judiciously chosen.  More
specifically, for $\lambda\in\Lambda^+$, let $g=g(\lambda)$ be the
solution of
\begin{equation}\label{e:g(lambda)}
(g,v)_\lambda
=
\int_{\Omega} 
\lambda (\nabla {g}\cdot \nabla v) + {g} v \, dx
=
e'(u;v),\quad\forall v\in H.
\end{equation}
\begin{remark}
\label{r:Riesz}
  In the following we will, \revt{in particular}, refer to the
  gradient \revtt{$g_1:=g(1)$} which corresponds to \revt{the usual inner
    product \eqref{e:(u,v)_Ht} and is also obtained by setting
    $\lambda=1$ in \eqref{e:g(lambda)}}, and to \revt{the gradient}
  \revtt{$g_0:=g(0)$} which corresponds to $\lambda=0$ \revt{in
    \eqref{e:g(lambda)}}.  Usually, \revt{$g_1$ and $g_0$ are referred
    to as, respectively, the $H^1$ and $L^2$ Riesz representations of
    $e'(u)$.}

In general $g_0\notin H^1_0(\Omega)$, but we have
\begin{eqnarray}
-\Delta g_1+g_1
&=&
g_0
\;\, \mbox{\it in ${\cal D'}(\Omega)$},
\label{e:g_1,g_0}\\
-\nabla\cdot(\lambda\nabla g(\lambda))+g(\lambda)
&=&g_0\;\, \mbox{\it in ${\cal D'}(\Omega)$}.
\label{e:g(l),g_0}
\end{eqnarray}
\end{remark}
Note that we will use the symbol $g$ to denote the operator
$\lambda\in\Lambda\mapsto g=g(\lambda)\in H$, or to denote an element
of $H$ --- the meaning will always be clear from the context.

Now assume we are at a certain iteration of the gradient method with $u$
known, which we seek to update to a new value $\widetilde{u}$, cf.~step
5 in Algorithm \ref{alg:gradient}. For this, first we look for a
certain $\hlambda \in \Lambda^+$, defined by\footnote{All along this
  paper, the symbol ``$\widehat{\phantom{x}}$'' will be used to denote
  the solution of a minimization problem, whereas the symbol
  ``$\widetilde{\phantom{x}}$'' will be used to represent an updated value
  of a variable.}
\begin{equation}\label{e:min(lambda)->1}
j(\hlambda) := \min\left\{
j(\lambda):=f\circ g(\lambda),\;\,  \lambda\in \Lambda^+
\right\}, \quad \text{\it where}
\quad
f(g):=\frac{\kappa}{2}e''(u;g,g)-e'(u;g).
\end{equation}
The reason for introducing the step size $\kappa$ in this equality
will be clear from Remark \ref{r:name} and also later during the error
analysis in Section \ref{sec:errors}.  Note that, if problem
\eqref{e:min(lambda)->1} has a solution $\hlambda \in \Lambda^+$, then
we will show (see Proposition \ref{p:g'->1}) that $\hlambda$ solves
\begin{equation}\label{e:e''(lambda)=e'}
e''(u;\kappa g(\hlambda),g'(\hlambda;\ell))=e'(u;g'(\hlambda;\ell)),\quad\forall \ell\in \Lambda,
\end{equation}
where $g'(\hlambda;\ell)$ denotes the derivative of $g$ at $\hlambda$
in the direction $\ell$.  Then, the modified gradient approach will
consist of Algorithm \ref{alg:gradient} with step 4 amended as follows
\begin{eqnarray}
&4.&\mbox{determine $g=g(\hlambda)$, where $\hlambda\in\Lambda^+$ is such that $g(\hlambda)\in H$ solves (\ref{e:e''(lambda)=e'})}.
\label{alg:gradient+2}
\end{eqnarray}
\begin{remark}\label{r:name}
  Clearly, our approach is equivalent to the gradient method, but with
  the classical inner product $(\cdot,\cdot)_H$ replaced with
  $(\cdot,\cdot)_\lambda$.

  From equation \eqref{e:e''(lambda)=e'} it follows that
  $e''(u;h-\kappa g,g'(\hlambda;\ell))=0$, where $h$ is
  Newton's step.  This means that $P_u(h-\kappa g)=0$,
  where $P_u:H\mapsto T_u$ is the projection from $H$ to $T_u$, in
  which $T_u= \Span\{g'(\hlambda;\ell),\; \ell\in\Lambda\}$ is the
  tangent space to the manifold $\{g(\lambda),\,
  \lambda\in\Lambda^+\}\subset H$ at $g(\hlambda)$, determined with
  respect to the inner product $e''(u;\cdot,\cdot)$.

  If $T_u=H$, then $\kappa g(\hlambda)=h$ and our method reduces to
  Newton's method.  However, here we have $dim(\Lambda)=N$, so that
  {in general} $T_u\neq H$ and $\kappa g$ will be close to $h$ in
  the sense that $P_u(h-\kappa g)=0$. This relation will be the key
  ingredient to prove {in the demonstration} that our gradient
  method, in addition to the linear convergence of a standard gradient
  method, has also a quadratic convergence in a certain sense
  depending on $T_u$ and the projection $P_u$.  This will be explained
  in the next sections.
%
\end{remark}

Our method critically depends on the choice of $\lambda$ and the
following proposition offers a first glimpse of what may happen with
the solution of problem \eqref{e:min(lambda)->1}.
\begin{proposition}\label{p:lambda^k->}
  Let $(\lambda^k)$ be a minimizing sequence of $j$ in $\Lambda^+$ and
  $(g(\lambda^k))$ be the \revt{corresponding} sequence of gradients
  \revt{defined in} \eqref{e:g(lambda)}.  Then, up to a subsequence,
  $(g(\lambda^k))$ converges weakly in $H^1_0(\Omega)$ and strongly
  in $L^2(\Omega)$ to \revt{an element} $g\in H$, while for the 
  sequence $(\lambda^k)$ one of the following cases may occur.\\
  (i) 
  There 
  exist $\hlambda=[\hlambda_1,\ldots,\hlambda_N]\in\Lambda^+$ and a subsequence of $(\lambda^k)$,
  still denoted $(\lambda^k)$, such that $\lim_{k\to\infty}\lambda^k_i=\hlambda_i$ for all
  $i=1,\dots,N$. 
  In this case $g=g(\hlambda)$, i.e. 
\begin{equation}\label{e:g(hlambda),1}
\int_\Omega \hlambda(\nabla g(\hlambda)\cdot\nabla v)+g(\hlambda)v \, dx = e'(u;v),\;\, \forall v\in H,
\end{equation}
and  $\hlambda$ solves \eqref{e:min(lambda)->1}.
\\
(ii) 
 There exist $\hlambda=[\hlambda_1,\ldots,\hlambda_N]\in\partial\Lambda^+$,
 $I_0\subset I$, $I_\infty\subset I$, with 
 $\hlambda_i=0$ for all $i\in I_0$,
 $\hlambda_i=+\infty$ for all $i\in I_\infty$,
 $0<\hlambda_i<+\infty$ for all $i\in I\backslash(I_0\cup I_\infty)$, 
 and a subsequence of $(\lambda^k)$, still denoted $(\lambda^k)$, such that 
 $\lim_{k\to\infty}\lambda^k_i=\hlambda_i$, for all $i=1,\dots,N$.
 In this case $g\in H^1_0(\Omega;\Omega_\infty^\text{const})$ solves
\begin{eqnarray}
\hspace*{-8mm}
\int_{\Omega\backslash \Omega_0}
\hlambda (\nabla g\cdot\nabla v) + gv  \,dx
+
\int_{\Omega_0} gv \, dx
&=& e'(u;v),
\quad
\forall v\in H^1_0(\Omega;\Omega_\infty^\text{const}),
\label{e:g-in-G-G0-Ginf}
\end{eqnarray}
where 
\begin{eqnarray}
\hspace*{-6mm}
\Omega_0&=&\cup\{\Omega_i,\; i\in I_0\},
\\
\hspace*{-6mm}
\Omega_\infty&=&\cup\{\Omega_i,\; i\in I_\infty\},
\\
\hspace*{-6mm}
H^1_0(\Omega;\Omega_\infty^\text{const})
&=&
\{
v\in H^1_0(\Omega),
\;\,
v=C_i\in\mathbb R\;\, in\;\, \Omega_i,\;\, \forall i\in I_\infty\}.
\label{e:H1_0(G;G0,Gi)}
\end{eqnarray}
Furthermore, if we define $g(\hlambda)=g$, with $g$ given by \eqref{e:g-in-G-G0-Ginf}, we have
\begin{equation}
j(\hlambda)\leq\liminf_{k\to\infty}j(\lambda^k).
\label{e:j(hl),=liminf}
\end{equation}

\end{proposition}
{\bf Proof}.
Let $(\lambda^k)$ be a minimizing sequence of $j$ in $\Lambda^+$ and 
$\lambda^k=[\lambda^k_1,\ldots,\lambda^k_N]$. 
Note that $g(\lambda^k)$ is well defined by 
\begin{equation}\label{e:g(lambda^k)}
\int_{\Omega}
\lambda^k(\nabla g(\lambda^k)\cdot\nabla v) + g(\lambda^k)v \, dx
= 
e'(u;v),\quad\forall v\in H.
\end{equation}
Note also that from the ellipticity of $f$ in $H$,
cf.~\eqref{e:e->C0}, \eqref{e:e->convex} and \eqref{e:min(lambda)->1},
it follows that the sequence $g(\lambda^k)$ is bounded in
$H^1(\Omega)$.  Therefore, up to a subsequence, we may assume that
$g(\lambda^k)$ converges weakly in $H$ and strongly in $L^2(\Omega)$
to a certain $g\in H$.

As ${\rm dim}(\Lambda)=N$, there exist $\hlambda=[\hlambda_1,\ldots,\hlambda_\revt{N}]$, with 
$\hlambda_i\in[0,+\infty]$, and a subsequence of $(\lambda^k)$, still denoted $(\lambda^k)$, 
such that $\lim_{k\to\infty}\lambda^k_i=\hlambda_i$ for all $i\in I$.
Two cases may occur.\\
1)
$\hlambda\in\Lambda^+$, i.e., $\hlambda_i\in(0,\infty)$ for all $i\in I$. 
From (\ref{e:g(lambda^k)}) we obtain
\begin{eqnarray*}
\int_{\Omega}
\hlambda (\nabla g\cdot\nabla v) + gv \,dx
&=& 
\int_{\Omega}(\hlambda-\lambda^k) (\nabla g\cdot\nabla v) \,dx
\\
&+&
\int_{\Omega}\lambda^k (\nabla(g-g(\lambda^k))\cdot\nabla v) + (g-g(\lambda^k))v \,dx
\\
&+&
\int_{\Omega}\lambda^k (\nabla g(\lambda^k)\cdot \nabla v) + g(\lambda^k)v \,dx
\\
&=&
\int_{\Omega}(\hlambda-\lambda^k) (\nabla g\cdot\nabla v) \,dx
\\
&+&
\int_{\Omega}\lambda^k (\nabla(g-g(\lambda^k))\cdot\nabla v) + (g-g(\lambda^k))v \,dx
\\
&+&
e'(u;v),\quad\forall v\in H.
\end{eqnarray*}
Then, letting $k$ go to infinity gives
\begin{eqnarray*}
\int_{\Omega}
\hlambda (\nabla g\cdot \nabla v) + gv \,dx
&=& 
e'(u;v),\quad\forall v\in H,
\end{eqnarray*}
which proves that $g=g(\hlambda)$.  Note that it is easy to show that 
the subsequence $(g(\lambda^k))$ converges to $g$ strongly in $H$, and therefore
$\hlambda$ is the solution of (\ref{e:min(lambda)->1}) because $j$ is continuous
in $H$.
%
\\
2) 
$\hlambda\in\partial\Lambda^+$.
\revt{Then,} there exist $I_0\subset I$ \revt{and} $I_\infty\subset I$
such that
$\hlambda_i=0$ for $i\in I_0$, 
$\hlambda_i=+\infty$ for $i\in I_\infty$, 
$\hlambda_i\in(0,+\infty)$ for $i\in I\backslash(I_0\cup I_\infty)$, 
and a subsequence of $(\lambda^k)$, still denoted by $(\lambda^k)$, such that
$\lim_{k\to\infty}\lambda^k_i=\hlambda_i$ for all $i\in I$.
From (\ref{e:g(lambda^k)}), for each $g(\lambda^k)$ and $v\in H^1_0(\Omega;\Omega_\infty^{const})$
we have
\begin{equation*}
\int_{\Omega\backslash \Omega_0}
\lambda^k (\nabla g(\lambda^k)\cdot \nabla v) + gv \, dx
+
\int_{\Omega_0}
\lambda^k (\nabla g(\lambda^k)\cdot \nabla v) +gv \, dx
= e'(u;v).
\end{equation*}
Then, letting $k$ go to infinity gives \eqref{e:g-in-G-G0-Ginf},
because $\nabla v = 0$ in $\Omega_\infty$.

To show that $g$ is constant on each $\Omega_i$, $i\in I_\infty$, we take
$v=g(\lambda^k)$ in \eqref{e:g(lambda^k)}, so that we obtain
\begin{eqnarray*}
\int_{\Omega_\infty}
\lambda^k |\nabla g(\lambda^k)|^2 + |g|^2 \, dx
&=&
e'(u;g(\lambda^k)) 
- \int_{\Omega\backslash\Omega_\infty} \lambda^k |\nabla g(\lambda^k)|^2 + |g|^2 \, dx
\\
&\leq&
C,
\end{eqnarray*}
because $(g(\lambda^k))$ is bounded in $H$ and $(\lambda^k)$ is
bounded in $L^\infty(\Omega\backslash\Omega_\infty)$.  It follows that
$\lim_{k\to\infty}|\nabla g(\lambda^k)|=0$ and
$\lim_{k\to\infty}g(\lambda^k)=g$ in $L^2(\Omega_\infty)$.  Hence,
$g=C_i$ in $\Omega_i$, $C_i\in\mathbb R$.  \revt{Thus, $g\in
  H^1_0(\Omega;\Omega_\infty^\text{const})$ solves
  \eqref{e:g-in-G-G0-Ginf}}.

Finally, \eqref{e:j(hl),=liminf} follows from the fact that $f$ is
convex and strongly continuous in $H$, so $f$ is weakly lower
semi-continuous, see \cite{brezis-1}.  \hfill$\Box$ 
\begin{remark}\label{r:I0,Iinf}
  While analyzing case (ii) we will use the following notation.  For
  $\lambda=[\lambda_1,\ldots,\lambda_N]\in\partial\Lambda^+$ we write
\begin{equation}
\left\{
\begin{array}{rcl}
 I_{0,\lambda}&=&\cup\{i\in I,\;\, \lambda_i=0\},\\
 \Omega_{0,\lambda}&=&\cup\{\Omega_i,\;\; i\in I_{0,\lambda}\},
\end{array}
\right.
\quad
\left\{
\begin{array}{rcl}
 I_{\infty,\lambda}&=&\cup\{i\in I,\;\, \lambda_i=+\infty\},\\
 \Omega_{\infty,\lambda}&=&\cup\{\Omega_i,\;\; i\in I_{\infty,\lambda}\}.
\end{array}
\right.
\end{equation}
For $\hlambda\in\partial\Lambda^+$ instead we write
\begin{equation}
\left\{
\begin{array}{rcl}
 I_0&=&\cup\{i\in I,\;\, \hlambda_i=0\},\\
 \Omega_0&=&\cup\{\Omega_i,\;\; i\in I_0\},
\end{array}
\right.
\quad
\left\{
\begin{array}{rcl}
 I_\infty&=&\cup\{i\in I,\;\, \hlambda_i=+\infty\},\\
 \Omega_\infty&=&\cup\{\Omega_i,\;\; i\in I_\infty\}.
\end{array}
\right.
\end{equation}
\end{remark}

\begin{remark}\label{r:hlambda,H10}
  In the case when $dim(\Lambda)=1$, i.e., $\Lambda=\mR$, and  the space $H$
 is equipped with the inner product
\begin{equation}
 (u,v)_H=\int_\Omega \nabla u \cdot\nabla v \, dx,
\label{eq:ipH10}
\end{equation}
the optimal weight $\hlambda$ is given explicitly. Indeed, if
$\displaystyle (u,v)_\lambda=\int_\Omega \lambda(\nabla u \cdot\nabla
v) \, dx$, then $g(\lambda)$ is defined by
\begin{equation*}
 \int_\Omega\lambda(\nabla g(\lambda)\cdot\nabla v) \, dx = e'(u;v).
\end{equation*}
This implies $g(\lambda)=\frac{1}{\lambda}g_1$ and then
\begin{equation*}
 j(\lambda)=\frac{\kappa}{2}\frac{1}{\lambda^2}e''(u;g_1,g_1)-\frac{1}{\lambda}e'(u;g_1).
\end{equation*}
It follows that the solution $\hlambda$ of \eqref{e:min(lambda)->1} is given by
\begin{equation}
\label{e:hlambda,N=1}
\hlambda=\kappa\frac{e''(u;g_1,g_1)}{e'(u;g_1)}. 
\end{equation}
Note that $\hlambda>0$ because $e'(u;g_1)>0$ and
$e''(u;g_1,g_1)>0$. Thus, in the case when $dim(\Lambda)=1$ and the
space $H$ is endowed with inner product \eqref{eq:ipH10}, the proposed
approach will consist of Algorithm \ref{alg:gradient} with step 4
amended as
\begin{eqnarray}
&4.&
\mbox{determine $g=g(\hlambda)$, where $\hlambda\in\Lambda^+$ is given by \eqref{e:hlambda,N=1}}.
\label{alg:gradient+2,N=1}
\end{eqnarray}
We remark that, interestingly, since $\hlambda$ is proportional to the
step size $\kappa$ and $g(\hlambda)$ is proportional to $1/\hlambda$,
in the present case the iterations produced by Algorithm
\ref{alg:gradient} will not depend on $\kappa$.

The optimal $\hlambda$ given in \eqref{e:hlambda,N=1} plays a similar
role to the parameter $\alpha$ used in the Barzilai-Borwein version of
the gradient method for minimization in $\mathbb \mR^n$
\cite{barzilai-1}. However, here the idea behind the choice of
$\hlambda$ given by \eqref{e:min(lambda)->1} or \eqref{e:hlambda,N=1}
is to approximate Newton's step. On the other hand, in
\cite{barzilai-1} the optimal $\alpha$ is chosen such that  the resulting gradient is a two-point
approximation to the secant direction used in the quasi-Newton methods.
\end{remark}

\section{Error analysis}
\label{sec:errors}
In the following, we first present the analysis of case (i) of
Proposition \ref{p:lambda^k->}.

\subsection{Error analysis: case $\hlambda\in\Lambda^+$}
\label{s:analysis->1}

\noindent
The following proposition gives the differentiability of the map $g$.
\begin{proposition}
\label{p:g'->1}
Let $\hlambda\in\Lambda^+$ be a solution of \eqref{e:min(lambda)->1}.
Then $g\in C^1(\Lambda,H)$ and $j\in C^1(\Lambda)$ near $\hlambda$.  
Furthermore, for all $v\in H$ and $\ell\in\Lambda$ we have
\begin{eqnarray}
\int_{\Omega}
\hlambda (\nabla g'(\hlambda;\ell)\cdot\nabla v) + g'(\hlambda;\ell)v \, dx 
&=&
- \int_{\Omega} \ell (\nabla g(\hlambda)\cdot \nabla v) \, dx,
\label{e:g'->1}\\
e''(u;\kappa g(\hlambda),g'(\hlambda;\ell))&=&e'(u;g'(\hlambda;\ell)),
\label{e:e''=e'->1}
\end{eqnarray}
where $g'(\hlambda;\ell)$ is the derivative of $g$ at $\hlambda$ in the direction $\ell$.
\end{proposition}
{\bf Proof}.
To prove the differentiability of $g$ we consider the map
\begin{eqnarray*}
F:\Lambda\times H&\mapsto&H'
\\
(\lambda,g)&\to& F(\lambda,g),
\quad 
F(\lambda,g)(v) 
=
\int_{\Omega}
\lambda (\nabla g\cdot\nabla v) + gv \, dx
-
e'(u;v),\;\, v\in H.
\end{eqnarray*}
Note that $F$ is $C^1$ and
\begin{equation*}
\partial_g F(\hlambda,g)(z)
= 
\int_{\Omega}
\hlambda (\nabla z\cdot\nabla v) + zv \, dx,\; z\in H.
\end{equation*}
It follows from the Lax-Milgram lemma that $\partial_g F(\hlambda,g)$
defines an isomorphism from $H$ to $H'$.  Then, the differentiability
of $g$ is easily deduced by using the implicit function theorem and
the fact that the \revt{equation} $F(\lambda,g)=0$ has a unique
solution $g\in H$ for any given $\lambda\in\Lambda^+$. In addition, it
follows that $\lambda\in\Lambda\mapsto j(\lambda)\in\mR$ is also $C^1$
near $\hlambda$, because $g\in C^1(\Lambda;H)$ and $f$ is continuous
in $H$.

Equalities \eqref{e:g'->1} and \eqref{e:e''=e'->1} are obtained after
straightforward computations.  
\hfill$\Box$

\begin{corollary}
\label{c:T->1}
Let $\hlambda\in\Lambda^+$ be a solution of \eqref{e:min(lambda)->1},
$T_u=\Span\{g'(\hlambda;\ell),\;\ell\in\Lambda\}$ and $P_u:H\mapsto T_u$
be the projection operator with respect to the inner product
$e''(u;\cdot,\cdot)$, i.e.,
\begin{equation}
e''(u;w-T_u w,v) = 0,\quad \forall w\in H,\;\, v\in T_u. 
\end{equation}
Then $P_uh = P_u(\kappa g)$ and $d\leq\dim(T_u)\leq {N}$, where
$d=\rank\{i,\; g(\hlambda)\neq g_0 \ \textrm{in} \ {\cal D}'(\Omega_i)
\}$ (we recall that $g_0$ is the $L^2$ \revt{representation of
  $e'(u)$, cf.~Remark \ref{r:Riesz}}).
\end{corollary} 
{\bf Proof}.  
From \eqref{e:e''=e'->1} and
$e''(u;h,v)=e'(u;v)$ for all $v\in H$, it follows $P_u h = P_u(\kappa
g)$.

Clearly $\dim(T_u)\leq N$.
Now we show that $d\leq \dim(T_u)$.  For simplicity and without loss
of generality we assume that $g(\hlambda)\neq g_0$ in ${\cal D}'(\Omega_i)$ for all $i=1,\ldots,d$.
%
It is enough to show that $\{g'(\hlambda;\ell^i),\; i=1,\ldots,d\}$
are linearly independent. Let 
$\sum_{i=1,d}\alpha_i g'(\hlambda;\ell^i)=0$, $\alpha_i\in\mR$. From \eqref{e:g'->1} we
obtain
\begin{eqnarray*}
0
&=&
\int_{\Omega}\hlambda 
\nabla \left(\sum_{i=1,d}\alpha_i g'(\hlambda;\ell^i)\right)\cdot\nabla v \, dx
+
\left(\sum_{i=1,d}\alpha_i g'(\hlambda;\ell^i)\right)v \, dx
\\
&=&
-
\int_{\Omega}
\left(\sum_{i=1,d}\alpha_i\ell^i\right)\nabla g(\hlambda)\cdot\nabla v \, dx.
\end{eqnarray*}
Then, taking $v\in{\cal D}(\Omega_i)$ gives
\[
0 
= 
\sum_{j=1}^d \alpha_j \int_{\Omega} \ell^j (\nabla g(\hlambda)\cdot\nabla v) \, dx 
=
\alpha_i \int_{\Omega_i} \nabla g(\hlambda)\cdot \nabla v \, dx
=
\alpha_i \int_{\Omega_i} \Delta g(\hlambda) v\, dx.
\]
Hence $\Delta g(\hlambda)=0$ in ${\cal D}'(\Omega_i)$.  Since
$-\nabla\cdot(\hlambda\nabla g(\hlambda)) +g(\hlambda)=g_0$ in ${\cal
  D}'(\Omega)$ and since $\hlambda$ is constant in each $\Omega_i$, it
follows \revt{that} $\alpha_i(g(\hlambda)-g_0)=0$ in ${\cal
  D}'(\Omega_i)$, hence $\alpha_i = 0$.  \hfill$\Box$

\begin{remark}
  The \revt{estimate of} the dimension of $T_u$ is optimal. In fact,
  we can prove that $d={\rm dim}\{i\in I,\;
  g'(\hlambda;\ell^i)\neq0\}$.  Indeed, if $\Delta g(\hlambda)=0$ in
  $\Omega_i$, we can \revt{show that} $\partial_\nu g(\hlambda)=0$ on
  $\partial\Omega_i$, \revt{where $\nu$ is the direction of the normal
    vector on $\partial\Omega_i$,} and then from (\ref{e:g'->1}) we
  get $g'(\hlambda;\ell^i)=0$.
\end{remark}

\noindent
Now we are able to prove the  error estimates for our method.
\begin{theorem}
\label{th:P+newton->1}
Assume $e(u)$ satisfies the assumptions of Theorems \ref{th:gradient}
and \ref{th:newton} near the solution $\hu$ of (\ref{e:min{e(u)}}).
Let $u$ be close to $\hu$ and $\widetilde{u}$ be given by Step 5 of
Algorithm \ref{alg:gradient} with $g=g(\hlambda)$ and
$\hlambda\in\Lambda^+$ a solution of (\ref{e:min(lambda)->1}). Then we
have
\begin{eqnarray}
\|\widetilde{u}-\hu\|_\hlambda
&\leq& 
\epsilon \|u-\hu\|_\hlambda,
\label{e:|tu-u|->1}\\
\|P_u(\widetilde{u}-\hu)\|_H
&\leq& 
C\|u-\hu\|_H^2,
\label{e:|P(tu-u)|->1}
\end{eqnarray}
with $\epsilon\in(0,1)$ depending only on $\kappa$ and $e$ and $C>0$
depending only on $e$.
\end{theorem} 
{\bf Proof}.  Estimate \eqref{e:|tu-u|->1} follows from
Theorem \ref{th:gradient}, where the norm is changed to
$\|\cdot\|_\hlambda$, because the gradient is now computed with respect to
the inner product $(\cdot,\cdot)_\hlambda$.  

For \eqref{e:|P(tu-u)|->1},
we note that $P_u$ is a linear continuous operator with 
$\|P_u\|_H\leq M$. Then
\begin{eqnarray*}
P_u(\widetilde{u}-\hu)
&=&
P_u(u-\kappa g(\hlambda) - \hu)
\\
&=&
P_u(u-\hu) - P_u(\kappa g(\hlambda)) \qquad(use\;\ Corollary\;\,\ref{c:T->1})
\\
&=&
P_u(u-\hu) - P_u(h)
\\&=&
P_u(u-h - \hu).
\end{eqnarray*}
Therefore
\begin{eqnarray*}
\|P_u(\widetilde{u}-\hu)\|_H
& = &
\|P_u(u-h-\hu)\|_H
\\
&\leq&
M\|u-h-\hu\|_H \qquad(use\; Theorem\; \ref{th:newton})
\\
&\leq&
C\|u-\hu\|_H^2.
\end{eqnarray*}
\begin{remark}
  \label{r:error->1}
  Theorem \ref{th:P+newton->1} states that $\|P_u(\widetilde{u}-\hu)\|_H$,
  the error of our method at a given step projected onto the tangent plane $T_u$,
  decreases at least quadratically in terms of $\|u-\hu\|_H$.
\end{remark}

\subsection{Error analysis: case $\hlambda\in\partial\Lambda^+$}
\label{s:analysis->2}

In case (ii) of Proposition \ref{p:lambda^k->} we are led to consider
$g(\hlambda)$ associated to $\hlambda\in\partial\Lambda^+$ with
$\hlambda_i=0$ for $i\in I_0$ and $\hlambda_i=+\infty$ for $i\in
I_\infty$, which solves \eqref{e:g(lambda)}.  To obtain error
estimates similar to the ones given by Theorem \ref{th:P+newton->1},
we would like to have differentiability results similar to the ones
given by Proposition \ref{p:g'->1}, which means that we \revt{would}
have to compare $g(\hlambda)$ with $g(\lambda)$,
$\lambda\in\partial\Lambda^+$.  \revt{However}, in general, for
$\lambda\in\partial\Lambda^+$ equation \eqref{e:g(lambda)} does not
provide \revt{an} estimate {in $H$} for $g(\lambda)$ and therefore
{the analysis from the previous section cannot be applied directly}.

\revt{On the other hand,} equation \eqref{e:g(lambda)} with
$\hlambda\in\partial\Lambda^+$ implies extra regularity for $e'(u)$,
in particular in $\Omega_0$. Assuming that $e''(u)$ {possesses} the
same kind of regularity, which comes naturally from the problem, we
will prove an error estimate for case (ii) of Proposition
\ref{p:lambda^k->} similar to the one already given in Theorem
\ref{th:P+newton->1}, but in a weaker norm.
\begin{proposition}\label{p:g=g_0}
Let 
$\Omega_{0,{H^1}}$ be the largest union of $\overline{\Omega}_i$ such that $g_0\in H^1(\Omega_{0,H^1})$
and $I_{0,H^1}=\{i\in I,\; \Omega_i\subset\Omega_{0,H^1}\}$.
\\
(i)
If $I_{0,H^1}=\emptyset$, then case (ii) of Proposition \ref{p:lambda^k->} does not happen.
\\
(ii)
If $I_{0,H^1}\neq\emptyset$, then $I_0\subset I_{0,H^1}$, $\Omega_0\subset \Omega_{0,H^1}$ and 
\begin{alignat}{2}
g(\hlambda)&=g_0 &  &\; \text{\it in} \ H^1(\Omega_0).
\label{e:g-in-G0}
\end{alignat}
(iii) Furthermore,
$e'(u;\cdot)$ is continuous in
$H^1_0(\Omega\backslash\Omega_0;\Omega_\infty^\text{const})\cap H^1(\Omega_0)$, where
\begin{eqnarray}
H^1_0(\Omega\backslash\Omega_0;\Omega_\infty^\text{const})\cap H^1(\Omega_0)
&=&
\{
v\in H^1(\Omega\backslash\Omega_0)\cap H^1(\Omega_0),\;\; 
\mbox{\it $v=0$\;\, on\;\, $\partial\Omega$\;\, and}	\nonumber\\
&&
\hspace*{2mm}
v=C_i\;\,\mbox{\it in}\; \Omega_\infty^i,\;\, i\in I_\infty\}.
\label{e:N1_0(G-G0;Ginf)}
\end{eqnarray}
\end{proposition} 
{\bf Proof}.  Indeed, from \eqref{e:g-in-G-G0-Ginf}
and \eqref{e:g(l),g_0}, we get $g(\hlambda)=g_0$ in ${\cal
  D}'(\Omega_0)$. As $g(\hlambda)\in H^1_0(\Omega)$, the claims (i) and
(ii) follow.

The form of the inner product $(g,v)_{\hlambda}$ and the fact that
$\hlambda=0$ in $\Omega_0$ \revt{imply} (iii).
\hfill$\Box$\\

\noindent
Motivated by Proposition \ref{p:lambda^k->}, we are led to the following definition.
\begin{definition}
  Let $\lambda=[\lambda_1,\ldots,\lambda_N]\in\partial\Lambda^+$ with 
  $I_{0,\lambda}\subset I_{0,H^1}$. 
  We define $g(\lambda)\in
  H^1_0(\Omega\backslash\Omega_{0,\lambda};\Omega_{\infty,\lambda}^\text{const})\cap H^1(\Omega_{0,\lambda})$ by
\begin{eqnarray}
\int_{\Omega\backslash\Omega_{0,\lambda}}
\lambda (\nabla g(\lambda)\cdot\nabla v) + g(\lambda)v \, dx
&+&
\int_{\Omega_{0,\lambda}}g(\lambda)v \, dx
=
e'(u;v),\nonumber\\
&&
\forall v\in H^1_0(\Omega\backslash\Omega_{0,\lambda};\Omega_{\infty,\lambda}^\text{const})\cap H^1(\Omega_{0,\lambda}).
\label{e:g''}
\end{eqnarray}
\end{definition}
\begin{proposition}
\label{p:g->2;wd}
Let  $\lambda\in\partial\Lambda^+$  with
with   $I_{0,\lambda}\subset I_{0,H^1}$. 
Then \eqref{e:g''} has a unique solution 
$g(\lambda)\in
H^1_0(\Omega\backslash\Omega_{0,\lambda};\Omega_{\infty,\lambda}^\text{const})\cap H^1(\Omega_{0,\lambda})$
and 
$g(\lambda)=g_0$ in $H^1(\Omega_{0,\lambda})$.
\end{proposition} 
{\bf Proof}.  
The existence and uniqueness of
$g(\lambda)$ follows from the Lax-Milgram lemma applied in the space 
$H^1_0(\Omega\backslash\Omega_{0,\lambda};\Omega_{\infty,\lambda}^\text{const})\cap L^2(\Omega_{0,\lambda})$
equipped with the inner product
\begin{equation*}
(g,v)_\lambda
=
\int_{\Omega}
\lambda (\nabla g\cdot\nabla v) + gv \, dx
=
\int_{\Omega\backslash\Omega_{0,\lambda}}
\lambda (\nabla u\cdot\nabla v) + gv \, dx
+
\int_{\Omega_{0,\lambda}}gv \, dx.
\end{equation*}
Reasoning as in \revt{case (ii) of} Proposition \ref{p:lambda^k->}, we
see that $g(\lambda)$ is constant in $\Omega_{\infty,\lambda}^i$, for
all $i\in I_{0,\lambda}$.  Finally, taking $v\in{\cal D}(\Omega_0)$ we
\revt{obtain} $\int_{\Omega_0}(g(\lambda)-g_0)v \,dx=0$, which implies
\revt{that} $g(\lambda)=g_0$ in ${\cal D}(\Omega_{0,\lambda})$ and, as
$g_0\in H^1(\Omega_{0,H^1})$, completes the proof.  \hfill$\Box$
\\

Returning to the minimization problem \eqref{e:min(lambda)->1} and in
view of \revt{case (ii) of} Proposition \ref{p:lambda^k->}, we are led to consider
the problem
\begin{equation}
\label{e:min(lambda)->2}
\text{find $\hlambda\in\partial\Lambda^+$ such that }\; 
j(\hlambda) := \min\{j(\lambda)=(f\circ g)(\lambda),\;\, \lambda\in\partial\Lambda^+\}
\end{equation}
and \revt{from it eventually obtain} a necessary condition
\revt{analogous to} \eqref{e:e''=e'->1}, which was a key element in
proving estimate \eqref{e:|P(tu-u)|->1}.

We would repeat the \revtt{analysis} already applied to problem
\eqref{e:min(lambda)->1}, \revtt{as} in Section \ref{s:analysis->1}.  
However,
since $g(\lambda)$ \revt{now} defined via \eqref{e:g''} does not in
general belong to $H^1_0(\Omega)$, $j(\lambda)$ may not be well
defined.

It appears that there are no general conditions on the data which
would ensure that $g(\lambda)\in H^1_0(\Omega)$ when
$\lambda\in\partial\Lambda^+$. We will thus proceed with the analysis
of this case under the following stronger assumptions on $e'$ and
$e''$, which are motivated by the continuity of $e'(u;\cdot)$ in
$H^1(\Omega\backslash\Omega_0)\cap H^1(\Omega_0)$, see Proposition
\ref{p:g=g_0}.

Let us \revt{introduce} the following definitions
\begin{eqnarray}
{\cal H}
&=&
\{u\in H^1(\Omega\backslash\Omega_{0,H^1}),\; v\in H^1(\Omega_i),\;\, i\in I_{0,H^1},\;\, 
u=0\;\; on\;\,\partial\Omega\},
\\
(u,v)_{{\cal H}}
&=&
\int_{\Omega\backslash\Omega_{0,H^1}}(\nabla u\cdot\nabla v) + uv \, dx
+
\sum_{i\in I_{0,H^1}}\int_{\Omega_i}(\nabla u\cdot\nabla v) + uv \, dx,
\\
\|v\|_{{\cal H}}^2
&=&
(v,v)_{{\cal H}}.
\end{eqnarray}
The set ${\cal H}$ equipped with the inner product
$(v,v)_{{\cal H}}$ is a Hilbert space.

In the reminder of this section we will assume
\begin{equation}\label{e:assumption(H10)}
\left\{
\begin{array}{l}
\mbox{\it 
$e$, $e'$ and $e''$ satisfy all the conditions of Theorems \ref{th:gradient} and \ref{th:newton} with}
\\
\mbox{\it 
$H$ replaced by ${\cal H}$}.
\end{array}
\right.
\end{equation}
Moreover, we will assume
\begin{eqnarray}
|e''(u;v,w)|&\leq&M\|v\|_{{\cal H}}\|w\|_{{\cal H}},
\label{e:e->C0,'}\\
e''(u;v,v)
&\geq&
m\|v\|_{{\cal H}}^2,
\label{e:e->convex,'}
\end{eqnarray}
with certain $0<m<M<\infty$.

\begin{proposition}
\label{p:min(lambda)->2}
Assume $e''$ satisfies \eqref{e:assumption(H10)}--\eqref{e:e->convex,'}. 
For $\lambda\in\partial\Lambda^+$ with $I_{0,\lambda}\subset I_{0,H^1}$ let
$g(\lambda)\in
H^1_0(\Omega\backslash\Omega_{0,\lambda};\Omega_{\infty,\lambda}^\text{const})\cap H^1(\Omega_{0,\lambda})$
be defined by (\ref{e:g''}).  
Then the problem \eqref{e:min(lambda)->2} has a solution $\hlambda\in\partial\Lambda^+$.
\end{proposition} 
{\bf Proof}. 
Let $(\lambda^k)$ be a sequence in
$\partial\Lambda^+$ minimizing $j$ in $\partial\Lambda^+$.  As
$\dim(\Lambda)=N$, without loss of generality, we may assume that
there exist $I_0\subset I$, $I_\infty\subset I$ such that
$I_{0,\lambda^k}=I_0$, $I_{\infty,\lambda^k}=I_\infty$ for all $k$.
It follows {that}
$\Omega_{0,\lambda^k}=\Omega_0$, $\Omega_{\infty,\lambda^k}=\Omega_\infty$.

Since $f$ is elliptic in ${\cal H}$ and $g(\lambda^k)\in
H^1(\Omega\backslash\Omega_0;\Omega_\infty^{\rm const})\cap
H^1(\Omega_0)$, for all $k$, necessarily $(g(\lambda^k))$ is bounded
in $H^1(\Omega\backslash\Omega_0)\cap H^1(\Omega_0)$. Therefore,
without loss of generality, we may assume that $(g(\lambda^k))$
converge weakly in $H^1(\Omega\backslash\Omega_0)\cap H^1(\Omega_0)$
and strongly in $L^2(\Omega)$ to a certain $g\in
H^1(\Omega\backslash\Omega_0)\cap H^1(\Omega_0)$.  As $f$ is convex,
it follows that $f(g)\leq\liminf_{k\to\infty}j(\lambda^k)$, \revt{see
  \cite{brezis-1}}.

To conclude that
\eqref{e:min(lambda)->2} has a solution, it is enough to show that
$g=g(\hlambda)$ for a certain $\hlambda\in\partial\Lambda^+$.
For the sequence $\lambda^k$ two cases may occur.\\
(i) 
There exist $\hlambda=[\hlambda_1,\ldots,\hlambda_N]$ with 
$\hlambda_i=0$ for $i\in I_0$,
$\hlambda_i=+\infty$ for $i\in I_\infty$,
$\hlambda_i\in(0,+\infty)$ for $i\in I\backslash(I_0\cup I_\infty)$,
and a subsequence of $(\lambda^k)$, still denoted $(\lambda^k)$,
such that 
$\lim_{k\to\infty}\lambda^k_i=\hlambda_i$ for all $i\in I$.
Note that $g_k=g(\lambda^k)\in
H^1(\Omega\backslash\Omega_0;\Omega_\infty^{\rm const})\cap H^1(\Omega_0)$ satisfies
\eqref{e:g''}, i.e.
\begin{eqnarray}
\int_{\Omega\backslash\Omega_0}
\lambda^k (\nabla g_k\cdot\nabla v) + g_kv \, dx
+
\int_{\Omega_0} g_kv \, dx
&=& e'(u;v),
\label{e:g'',k}
\end{eqnarray}
for all $v\in H^1(\Omega\backslash\Omega_0;\Omega_\infty^{\rm const})\cap H^1(\Omega_0)$.
Passing {to the} limit in 
\eqref{e:g'',k},  we find that $g\in H^1_0(\Omega\backslash\Omega_0;\Omega_\infty^\text{const})\cap H^1(\Omega_0)$ 
solves \eqref{e:g''} so $g=g(\hlambda)$.
\\
(ii)
There exist $\hlambda=[\hlambda_1,\ldots,\hlambda_N]$, $i_0\subset I$, $i_\infty\subset I$ with 
$\hlambda_i=0$ for $i\in I_0\cup i_0$,
$\hlambda_i=+\infty$ for $i\in I_\infty\cup i_\infty$,
$\hlambda_i\in(0,+\infty)$ for $i\in I\backslash((I_0\cup i_0)\cup (I_\infty\cup i_\infty))$,
and a subsequence of $(\lambda^k)$, still denoted $(\lambda^k)$,
such that 
$\lim_{k\to\infty}\lambda^k_i=\hlambda_i$ for all $i\in I$.

We take $v\in H^1_0(\Omega\backslash\Omega_0;(\Omega_\infty\cup\omega_\infty)^\text{const})\cap H^1(\Omega_0)$
in  \eqref{e:g'',k}, where 
$\omega_0=\cup\{\Omega_i,\; i\in i_0\}$,
$\omega_\infty=\cup\{\Omega_i,\; i\in i_\infty\}$,
and we obtain
\begin{eqnarray*}
\int_{\Omega\backslash(\Omega_0\cup \omega_0)}
\lambda^k (\nabla g_k\cdot\nabla v) + g_kv \, dx
&+&
\int_{\Omega_0}
g_kv \, dx
+
\int_{\omega_0}
\lambda^k (\nabla g_k\cdot\nabla v) + g_kv \, dx
\nonumber\\
&=& 
e'(u;v).
\label{e:g'',k'}
\end{eqnarray*}
Letting $k\to\infty$ gives
\begin{eqnarray*}
\int_{\Omega\backslash(\Omega_0\cup \omega_0)}
\hspace*{-5mm}
\hlambda (\nabla g\cdot\nabla v) + gv \, dx
+
\int_{\Omega_0\cup\omega_0}
\hspace*{-3mm}
gv \, dx
&=& 
e'(u;v),
\end{eqnarray*}
for all $v\in
H^1_0(\Omega\backslash\Omega_0;(\Omega_\infty\cup\omega_\infty)^\text{const})\cap
H^1(\Omega_0)$.  Reasoning as in Proposition \ref{p:lambda^k->}, we
find that $g=g(\hlambda)\in
H^1_0(\Omega\backslash\Omega_0;(\Omega_\infty\cup\omega_\infty)^\text{const})\cap
H^1(\Omega_0)$ solves \eqref{e:g''} with $\lambda=\hlambda$ and
$\Omega_0\cup\omega_0$ (respectively,
$\Omega_\infty\cup\omega_\infty$) instead of $\Omega_0$ (respectively,
$\Omega_\infty$).  \hfill$\Box$

\begin{remark}\label{r:sigma(lambda)}
  For $\hlambda\in\partial\Lambda^+$, in order to control the
  variations of $\hlambda$ in the set
  $\Omega\backslash(\Omega_0\cup\Omega_\infty)$ we consider ${\mathbb
    1}_{0,\infty}\in\partial\Lambda^+$ defined by
\[
{\mathbb 1}_{0,\infty}(x) 
=
\left\{
\begin{array}{ll}
0,&x\in \Omega_0\cup\Omega_\infty,\\
1,&x\in\Omega\backslash(\Omega_0\cup\Omega_\infty).
\end{array}
\right.
\]
Then we perturb $\hlambda$ with the elements of 
$\revt{\widehat{\Lambda}}:={\mathbb 1}_{0,\infty}\cdot \Lambda
=
\{{\mathbb 1}_{0,\infty}\cdot\lambda,\; \lambda\in\Lambda\}$.
\end{remark}
\begin{proposition}\label{p:g'->2}
Let $\hlambda\in\partial\Lambda^+$ be the solution of \eqref{e:min(lambda)->2} as given by Proposition 
\ref{p:min(lambda)->2}.
The map 
$\lambda\in\hat{\Lambda}
\mapsto 
g(\lambda)\in H^1_0(\Omega\backslash\Omega_0;\Omega_\infty^{\rm const})\cap H^1(\Omega_0)$
is $C^1$ near $\hlambda$.
Furthermore, $g'(\hlambda;\ell)\in H^1_0(\Omega\backslash\Omega_0;\Omega_\infty^{\rm const})\cap H^1(\Omega_0)$, 
where $g'(\hlambda;\ell)$ is the derivative of $g(\hlambda)$ at $\hlambda$ in \revt{the} 
direction $\ell\in\hat{\Lambda}$, and satisfies 
\begin{eqnarray}
\int_{\Omega\backslash\Omega_0}
\hlambda(\nabla g'(\hlambda;\ell)\cdot\nabla v)
&+&
g'(\hlambda;\ell) v \, dx
+
\int_{\Omega_0}
g'(\hlambda;\ell) v \, dx
=
-\int_{\Omega\backslash\Omega_0} 
\ell (\nabla g(\hlambda)\cdot\nabla v) \, dx,
\nonumber\\
&&
\forall
v\in H^1_0(\Omega\backslash\Omega_0;\Omega_\infty^{\rm const})\cap H^1(\Omega_0).
\label{e:Dg''}
\end{eqnarray}
In particular, $g'(\hlambda;\ell)=0$ in $\Omega_0$ and  for every $\ell\in \revt{\widehat{\Lambda}}$ we have
\begin{equation}\label{e:e''=e'->2}
 e''(u;\kappa g(\hlambda),g'(\hlambda;\ell))=e'(u;g'(\hlambda;\ell)).
\end{equation}
\end{proposition} 
{\bf Proof}.  The differentiability of $g$ is
deduced from the implicit mapping theorem as follows. Consider the map
$F$
\begin{equation*}
 \begin{array}{lcll}
  F:&\hat{\Lambda}\times H^1_0(\Omega\backslash\Omega_0;\Omega_\infty^\text{const})\cap L^2(\Omega_0)
  &\mapsto &
  (H^1_0(\Omega\backslash\Omega_0;\Omega_\infty^\text{const})\cap L^2(\Omega_0))'
  \\
  &(\lambda,g)&\to&
F(\lambda,g),
 \end{array}
\end{equation*}
with
\begin{equation*}
F(\lambda,g) 
= 
\int_{\Omega\backslash\Omega_0}
\lambda(\nabla g(\lambda)\cdot\nabla v)
+
g(\lambda) v \, dx
+
\int_{\Omega_0}
g(\lambda) v \, dx
-
e'(u;v).
\end{equation*}
Clearly, $F$ is $C^1$ near $(\hlambda,g(\hlambda))$. Furthermore, we have
\begin{eqnarray*}
\partial_g F(\hlambda,g(\hlambda))(z)
&=&
\int_{\Omega\backslash\Omega_0}
\hlambda(\nabla z\cdot\nabla v)
+
z v \, dx
+
\int_{\Omega_0}
zv \, dx,
\end{eqnarray*}
which defines an isomorphism from
$H^1_0(\Omega\backslash\Omega_0;\Omega_\infty^\text{const})\cap
L^2(\Omega_0)$ to its dual.  Then, the implicit mapping theorem and
the fact that $F(\lambda,g)=0$ has a unique solution $g\in
H^1(\Omega\backslash\Omega_0;\Omega_\infty^\text{const})\cap
H^1(\Omega_0)$ for any $\lambda\in \revt{\widehat{\Lambda}}$ gives the
differentiability of the map $g$.  Note that, a priori, \revt{the}
implicit mapping theorem \revt{ensures} the differentiability of the
map $g$ in $L^2(\Omega_0)$. \revt{Then,} as $g(\lambda)=g_0\in
H^1(\Omega_0)$, see Proposition \ref{p:g->2;wd}, the differentiability
of the map $g$ in $H^1(\Omega_0)$ \revt{follows as well}.

Next, by direct computations we can easily show \eqref{e:Dg''}.
Furthermore, $g'(\hlambda;\ell)=0$ in $\Omega_0$ because
$g(\lambda)=g_0$ in $H^1(\Omega_0)$.

In regard to \eqref{e:e''=e'->2}, we recall that $e'(u)$ and $e''(u)$
are, respectively, linear and bilinear, and continuous in ${\cal H}$,
which together with the identity $g(\lambda)=g_0$ in $\Omega_0$ and
the differentiability of $g(\lambda)$ \revt{implies} the
differentiability of $\lambda\mapsto e'(u;g(\lambda))$ and
$\lambda\mapsto e''(u;g(\lambda),g(\lambda))$.  Then,
(\ref{e:e''=e'->2}) follows from straightforward computations.
\hfill$\Box$
\\

The error estimates are obtained in an analogous way to the
corresponding results in Section \ref{s:analysis->1}. First, we have a
result similar to Corollary \ref{c:T->1}.

\begin{corollary}\label{c:T->2}
Assume $e''(u)$ satisfies conditions \eqref{e:assumption(H10)}--\eqref{e:e->convex,'}.  
Let $\hlambda\in\partial\Lambda^+$ be a solution of
\eqref{e:min(lambda)->2}, $g=g(\hlambda)$, $T_u=\Span\{g'(\hlambda;\ell),\;\ell\in\hat{\Lambda}\}$ and
$P_u:{\cal H}\mapsto T_u$ be
the projection operator with respect to the inner product
$e''(u;\cdot,\cdot)$, i.e.,
\begin{equation}
e''(u;w-P_uw,v) = 0,\quad \forall w\in {\cal H},\;\, \forall v\in T_u.
\end{equation}
Then, $P_u h = P_u(\kappa g)$ and $d\leq \dim(T_u)\leq
N-|I_0|-|I_\infty|$, where $d=\rank\{i,\; g\neq g_0\; \mbox{in}\;
{\cal D}'(\Omega_i)\}$.

\end{corollary}
{\bf Proof}.
From \eqref{e:e''=e'->2} and the relation $e''(u;h,v)=e'(u;v)$ for all $v\in H$, 
it follows that $P_u h = P_u(\kappa g)$. 

Clearly $\dim(T_u)\leq N-|I_0|-|I_\infty|$.
Now we show that $\dim(T_u)\geq d$.  We assume that $g(\hlambda)\neq g_0$ in ${\cal D}'(\Omega_i)$, 
for all $i=1,\ldots,d$.  Let $\ell^i\in\hat{\Lambda}$,
$\ell^i(\Omega_j)=\delta_{ij}$.
It is enough to demonstrate that $\{g'(\hlambda;\ell^i),\; i=1,\ldots,d\}$
are linearly independent.  Let $\sum_{i=1,d}\alpha_i
g'(\hlambda;\ell_i)=0$, $\alpha_i\in\mR$. Then,
\begin{eqnarray*}
0
&=&
\int_{\Omega\backslash\Omega_0}
\hlambda 
 \nabla \left(\sum_{i=1,d}\alpha g'(\hlambda;\ell^i)\right)\cdot\nabla v 
+
\left(\sum_{i=1,d}\alpha_i \nabla g'(\hlambda;\ell^i)\right)v \, dx
\\
&+&
\int_{\Omega_0}
\left(\sum_{i=1,d}\alpha_i \nabla g'(\hlambda;\ell^i)\right)v \, dx
\\
&=&
\int_{\Omega\backslash\Omega_0}
\left(\sum_{i=1,d}\alpha_i \ell^i\right) (\nabla g(\hlambda)\cdot\nabla v) \, dx.
\end{eqnarray*}
In the equality above we take $v\in{\cal D}(\Omega_i)$. Then
\begin{eqnarray*}
0
&=&
\int_{\Omega\backslash\Omega_0}
\alpha_i(\nabla g(\hlambda)\cdot \nabla v) \, dx
=
- \alpha_i\int_{\Omega}\Delta g(\hlambda) v \, dx
=
- \alpha_i\int_{\Omega}(g(\hlambda)-g_0)v \, dx,
\end{eqnarray*}
because $-\hlambda\Delta g(\hlambda)+g(\hlambda)=g_0$ in $\Omega_i$,
which implies $\alpha_i=0$ and proves the claim.
\hfill$\Box$
\\

\noindent
Finally, we are able to prove the error estimate for the case
$\hlambda\in\partial\Lambda^+$.
\begin{theorem}
\label{th:P+newton->2}
Assume $e$ satisfies the conditions of Theorems \ref{th:gradient},
\ref{th:newton} with ${\cal H}$ instead of $H$, and $e''(u)$ satisfies
conditions \eqref{e:assumption(H10)}--\eqref{e:e->convex,'}.  Let
$\hlambda\in\partial\Lambda^+$ be a solution of
\eqref{e:min(lambda)->2}, as given by Proposition
\ref{p:min(lambda)->2}.  If $\widetilde{u}$ and $u$ are given as in
Algorithm \ref{alg:gradient} with $g=g(\hlambda)$, then
\begin{eqnarray}
\|\widetilde{u}-\hu)\|_\hlambda
&\leq& 
\epsilon \|u-\hu\|_\hlambda,
\label{e:|tu-u|->2}\\
\|P_u(\widetilde{u}-\hu)\|_{{\cal H}}
&\leq& 
C\|u-\hu\|_{{\cal H}}^2,
\label{e:|P(tu-u)|->2}
\end{eqnarray}
with 
$\epsilon\in(0,1)$  depending only on $\kappa$ and $e$ and $C>0$ depending only on $e$.
\end{theorem} 
{\bf Proof}.  The proof is analogous to the proof of
Theorem \ref{th:P+newton->1}.  
\hfill$\Box$

\begin{remark}
\label{r:error->2}
The estimates in Theorem \ref{th:P+newton->2} are similar to the ones
in Theorem \ref{th:P+newton->1}. However, in general, the quadratic
convergence established in Theorem \ref{th:P+newton->2} is slower than
the one provided by Theorem \ref{th:P+newton->1}, because in Theorem
\ref{th:P+newton->2} the dimension of the space $T_u$ is in general
smaller than the dimension of the space $T_u$ in Theorem
\ref{th:P+newton->1}.  Finally, estimate \eqref{e:|P(tu-u)|->2} {might
  be of no} interest because we may have
$\Omega\backslash(\Omega_0\cup\Omega_\infty)=\emptyset$ and so
$T_u=\{0\}$. In \revt{such} case $g=g_0$ \revt{would be} the $L^2$
gradient \revt{defined in the entire domain} $\Omega$ and the
question of its impact on the \revt{performance} of the gradient
method is open.

\revt{Lastly}, Theorem \ref{th:P+newton->2} \revt{provides an estimate
  applicable at a single step of the gradient approach, cf.~Algorithm
  \ref{alg:gradient},} where a certain $u$ is given and the regularity
of $e$ is \revt{determined} in terms of the set where the $L^2$
gradient $g_0=g_0(u)$ is $H^1$.  In order to \revt{be able to apply}
Theorem \ref{th:P+newton->2} at each step, \revt{one} should rather
consider \revt{iterations in the space} ${\cal H}=\{u\in
H^1(\Omega_i),\; i\in I,\; u=0\; on\; \partial\Omega\}$ and impose the
same assumptions as in Theorem \ref{th:P+newton->2}.

\end{remark}

\section{Determination of optimal weights $\hlambda$ and the corresponding gradients}
\label{sec:optimal}

In this section we describe the computational approach which can be
used to determine the optimal form of the inner product
\eqref{e:(u,v)_lambda}, encoded in its weight $\lambda$, and the
corresponding gradient $g(\lambda)$ at each iteration, cf.~modified
step 4 of Algorithm \ref{alg:gradient} given by
\eqref{alg:gradient+2}. We will focus on the case when $\hlambda \in
\Lambda^+$, cf.~Section \ref{s:analysis->1}, and in order to ensure
non-negativity of the weight, in our approach we will use the
representation $\lambda(x) = \eta^2(x)$, $\forall x \in \Omega$, where
$\eta \; : \; \Omega \mapsto \mR$ is a function {defined} below.
For consistency with the notation introduced in the previous sections
and without risk of confusion, hereafter we will use both $\lambda$
and $\eta$. Relation \eqref{e:g(lambda)} can then be expressed in the
strong form as
\begin{equation}
\left\{
\begin{alignedat}{2}
g -\nabla\cdot(\eta^2 \nabla g) &= g_0(u) & \quad & \text{\it in} \ \Omega, \\
g & = 0 & \quad & \text{\it on} \ \partial\Omega,
\end{alignedat}
\right.
\label{eq:g(eta)}
\end{equation}
where here $g=g(\eta)$, whereas the minimization problem \eqref{e:min(lambda)->1} becomes
\begin{equation}
\label{e:min(eta)->1}
j(\heta) := \min\left\{
j(\eta):=f\circ g(\eta),\;\,  \eta^2 \in \Lambda^+
\right\}.
\end{equation}
We will assume that the function $\eta(x)$ is represented with the ansatz
\begin{equation}
\eta(x) = \sum_{i=1}^N \eta_i \, \ell^i(x),
\label{eq:etaN}
\end{equation}
where $\{ \ell^i\}_{i=1}^N$ is a set of suitable basis functions.
Since in a fixed basis the function $\eta(x)$ is determined by the real
coefficients $\{ \eta_i\}_{i=1}^N$, we will also use the notation $\eta = [
\eta_1,\dots,\eta_N ]$. We thus obtain a finite-dimensional
minimization problem
\begin{equation}
\min\{ j(\eta),\;\, \eta=[ \eta_1,\dots,\eta_N ] \in \mR^N \}.
\label{eq:minj(eta)}
\end{equation}
Its minimizers $\heta = [ \heta_1,\dots,\heta_N ]$ satisfy the
following optimality conditions, which can be viewed as a discrete form
of \eqref{e:e''(lambda)=e'},
\begin{equation}
\left[ F_i(\heta) \right]
:
= 
\left[\Dpartial{j}{\eta_i}(\heta)\right]
=
[e''(u;\kappa g(\heta),g'_i(\heta))-e'(u;g'_i(\heta))]
=
[0], \qquad i=1,\dots,N,
\label{eq:alpha->F=0}
\end{equation}
where $F=[F_1,\ldots,F_N] \; : \; \mR^N \rightarrow \mR^N$ and 
$g'_i =g'_i(\eta)=
\left[\Dpartial{g}{\eta_i}(\heta)\right]$ satisfy the equations
\begin{equation}
\left\{
\begin{alignedat}{2}
g'_i - \nabla\cdot\left(\heta^2 \, \nabla g'_i \right) 
&= 2 \nabla\cdot\left( \heta \, \ell_i \nabla g\right) & \qquad & \mbox{\it in } \Omega, \\
g'_i &= 0 & & \mbox{\it on }\partial\Omega.
\end{alignedat}
\right.
\label{eq:g'i}
\end{equation}
The optimal weight $\heta$ can be found either by
directly minimizing $j(\eta)$, cf.~\eqref{eq:minj(eta)}, using a
version of the gradient-descent method, or by solving the optimality
conditions \eqref{eq:alpha->F=0} with a version of Newton's method.
The two approaches are described below.

\subsection{Optimal weights via gradient minimization}
\label{sec:etagrad}

While in practice one may prefer to use a more efficient minimization
approach, such as, e.g., the nonlinear conjugate-gradients method
\cite{nw00}, for simplicity of presentation here we focus on the
gradient steepest-descent method.  The step size $\tau_k$ along the
gradient can be determined by solving a line-minimization problem,
which can be done efficiently using for example Brent's method
\cite{nw00}. Step 4 of Algorithm \ref{alg:gradient},
cf.~\eqref{alg:gradient+2}, is then realized by the operations
summarized as Algorithm \ref{alg:etagrad}.  In actual computations it
may also be beneficial to prevent any of the values $\eta_i$ from
becoming too close to zero, which is achieved easily by imposing a
suitable bound on the step size $\tau_n$.
Having in mind the complexity analysis presented in Section
\ref{sec:complexity}, the termination condition for the main loop in
Algorithm \ref{alg:etagrad} is expressed in terms of the maximum
number $N_g$ of iterations, although in practice it will be more
convenient to base this condition on the relative decrease of
$j(\eta)$.
\begin{algorithm}[h!]
\begin{algorithmic}[1]
\STATE evaluate adjoint states $z$ and $z(g)$ (if $e(u)$ depends on
a PDE equation)
\STATE set $k = 1$
\REPEAT 
\STATE  evaluate $g(\eta)$ by solving \eqref{eq:g(eta)}
\STATE  evaluate $g'_i(\eta)$, $1=1,\dots,N$, by solving problems \eqref{eq:g'i}
\STATE  evaluate $e'(u; g'_i(\eta))$, $i=1,\ldots,N$
\STATE  evaluate $e''(u; \kappa g(\eta), g'_i(\eta))$, $i=1,\ldots,N$
\STATE  evaluate $F(\eta)$, cf.~\eqref{eq:alpha->F=0}
\STATE  perform line-minimization to determine optimal step size \newline
    \hspace*{1.0cm} $\tau_k = \argmin_{\tau>0} j(\eta - \tau F(\eta))$  \qquad (Brent's method \cite{nw00})
    \\
\STATE  set $\heta = \eta - \tau_k F(\eta)$
\STATE  set $\eta = \heta$
\STATE  set $k = k + 1$
\UNTIL{ \  $k = N_g$} 
\STATE obtain $g(\heta)$ by solving \eqref{eq:g(eta)} with $\eta = \heta$
\end{algorithmic}
\caption{Determination of optimal weight $\hlambda$ via gradient minimization \newline
     \textbf{Input:} \newline
 \hspace*{0.5cm} $N$ --- dimension of the space in which optimal weights are sought \newline
 \hspace*{0.5cm} $u \in H$ --- current approximation of minimizer $\hu$ \newline
 \hspace*{0.5cm} $\kappa > 0$ --- step size in the outer loop (Algorithm \ref{alg:gradient}) \newline
 \hspace*{0.5cm} $\{ \ell^i\}_{i=1}^N$ --- basis function for ansatz \eqref{eq:etaN} \newline
 \hspace*{0.5cm} $N_g$ --- maximum number of gradient iterations \newline
 \hspace*{0.5cm} $\eta$ --- initial guess for the weight \newline
 \textbf{Output:} \newline
 \hspace*{0.5cm} $\heta$ --- optimal weight \newline
 \hspace*{0.5cm} $g(\heta)$ --- corresponding optimal gradient
}
\label{alg:etagrad}
\end{algorithm}

\subsection{Optimal weights via Newton's method}
\label{sec:etanewton}

In addition to the gradient of $j(\eta)$ already given in
\eqref{eq:alpha->F=0}--\eqref{eq:g'i}, the key additional step
required for Newton's method is the evaluation of the Hessian of
$j(\eta)$, i.e.,
\begin{align}
\left[ \partial_jF_i(\eta) \right]
& = \left[\Dpartialmix{j}{\eta_i}{\eta_j}(\eta)\right] \qquad\qquad \qquad\qquad 
(i,j=1,\dots,N) \nonumber \\
& = \kappa
(e''(u;g'_j(\eta),g'_i(\eta)) + e''(u;g(\eta),g''_{ij}(\eta)))
- 
e'(u;g''_{ij}(\eta)), 
\label{eq:DF}
\end{align}
where $g(\eta)$ is given by \eqref{eq:g(eta)}, $g'_i(\eta)$ is given by \eqref{eq:g'i} and
$g''_{ij}=g''_{ij}(\eta) = \left[\Dpartialmix{g}{\eta_i}{\eta_j}(\heta)\right]$ satisfy the equations
\begin{equation}
\begin{alignedat}{2}
g''_{ij}(\eta)-\nabla\cdot(\eta^2\nabla g''_{ij}(\eta))
&= 2\left( \nabla\cdot(\ell_j \eta\nabla g'_i(\eta)) +  \nabla\cdot(\ell_i \eta\nabla g'_j(\eta))\right)  \\
&+ 2 \nabla\cdot(\ell_i \ell_j\nabla g(\eta)) &
\quad & \text{\it in} \ \Omega,\\
g''_{ij}(\eta)
& = 0 &
\quad & \text{\it on} \ \partial\Omega. 
\end{alignedat}
\label{eq:g''}
\end{equation}
For brevity, Newton's approach is stated in Algorithm
\ref{alg:etanewton} in its simplest form and in practice one would
typically use its damped (globalized) version in which the step along
Newton's direction $-\left[ DF(\eta) \right]^{-1}\cdot F(\eta)$ may be
reduced to ensure the residual $\| F(\eta)\|_2$ of equation
\eqref{eq:alpha->F=0} decreases between iterations \cite{k03}. A
similar step-size limitation may also be imposed in order to prevent
any of the values $\eta_i$ from becoming too close to zero. In
addition, in practice, a termination criterion based on the residual
$\| F(\eta)\|_2$ will be more useful. The criterion involving the
total number of iterations $N_n$ is used in Algorithm
\ref{alg:etanewton} only to simplify the complexity analysis which is
presented next.
\begin{algorithm}[h!]
\begin{algorithmic}[1]
\STATE evaluate adjoint states $z$ and $z(g'_i)$, $i=1,\ldots,N$ (if $e(u)$ depends on
a PDE equation)
\STATE set $k = 1$
\REPEAT 
\STATE  evaluate $g(\eta)$ by solving \eqref{eq:g(eta)}
\STATE  evaluate $g'_i(\eta)$, $i=1,\dots,N$, by solving \eqref{eq:g'i}
\STATE  evaluate $g''_{ij}(\eta)$, $i,j=1,\dots,N$, by solving \eqref{eq:g''}
\STATE  evaluate $e'(u; g'_i(\eta))$, $i=1,\ldots,N$
\STATE  evaluate $e'(u; g'_{ij}(\eta))$, $i,j=1,\ldots,N$
\STATE  evaluate $e''(u; g'_i(\eta),g'_j(\eta))$, $i,j=1,\ldots,N$
\STATE  evaluate $e''(u; g(\eta),g'_{ij}(\eta))$, $i,j=1,\ldots,N$
\STATE  evaluate the function $F(\eta)$, cf.~\eqref{eq:alpha->F=0}
\STATE  evaluate the Hessian $DF(\eta)$, cf.~\eqref{eq:DF}
\STATE  set $\heta = \eta - \left[ DF(\eta) \right]^{-1}\cdot F(\eta)$
\STATE  set $\eta = \heta$
\STATE  set $k = k + 1$
\UNTIL{ \  $k = N_n$} 
\STATE obtain $g(\heta)$ by solving \eqref{eq:g(eta)} with $\eta = \heta$
\end{algorithmic}
\caption{Determination of optimal weight $\hlambda$ using Newton's method \newline
     \textbf{Input:} \newline
 \hspace*{0.5cm} $N$ --- dimension of the space in which optimal weights are sought \newline
 \hspace*{0.5cm} $u \in H$ --- current approximation of minimizer $\hu$ \newline
 \hspace*{0.5cm} $\kappa > 0$ --- step size in the outer loop (Algorithm \ref{alg:gradient}) \newline
 \hspace*{0.5cm} $\{ \ell^i\}_{i=1}^N$ --- basis function for ansatz \eqref{eq:etaN} \newline
 \hspace*{0.5cm} $N_n$ --- maximum number of Newton iterations \newline
 \hspace*{0.5cm} $\eta$ --- initial guess for the weight \newline
 \textbf{Output:} \newline
 \hspace*{0.5cm} $\heta$ --- optimal weight \newline
 \hspace*{0.5cm} $g(\heta)$ --- corresponding optimal gradient
}
\label{alg:etanewton}
\end{algorithm}

\FloatBarrier

\subsection{Complexity analysis}
\label{sec:complexity}

In this section we estimate the computational cost of a single
iteration of Algorithms \ref{alg:etagrad} and \ref{alg:etanewton} in
which the optimal weight $\hlambda$ is computed using gradient
minimization and Newton's method, respectively, as described in
Sections \ref{sec:etagrad} and \ref{sec:etanewton}.  This cost will be
expressed in terms of:
(i)
the number $N$ of the degrees of freedom
characterizing the dimension of the weight space $\Lambda$,
cf.~\eqref{eq:etaN};
(ii)
the number $M$ {determining the cost of} the numerical solution of 
the elliptic boundary-value problems \eqref{eq:g(eta)}, \eqref{eq:g'i},
\eqref{eq:g''}; this latter quantity can be viewed as the number of
computational elements used to discretize the domain $\Omega$ (such as
finite elements/volumes, grid points or spectral basis functions);
(iii)
the number $K$ which is the typical number of line-search
iterations (line 9 in Algorithm \ref{alg:etagrad}). In the following we
will assume that the constants $C_1, C_2, \dots$ are all positive and
$\cO(1)$.

Both algorithms require first the evaluation of $e'(u;\ell_i)$, $i=1,\ldots,N$. 
{In general, the linear form can be expressed as}
\begin{equation}
e'(u;v) = \int_{\Omega} z v \, dz
\label{eq:z}
\end{equation}
{and, assuming that $v$ is already available, the cost of its
  approximation is determined by the cost of evaluating $z$ on
  $\Omega$ and the cost of the quadrature which is typically
  $\cO(M)$.}  If $z$ is a function given explicitly in terms of $u$,
then {it can be evaluated on $\Omega$ in terms of $\cO(M)$ operations.
  However,} in general, {when the energy depends on $u$ through
  some PDE,} which is the case of interest {here,} $z$ will be
given in terms of the solution of a suitably-defined adjoint PDE
problem. {Then,} for example, if the governing system is an
elliptic PDE problem in dimension $(d+1)$ with $u$ acting as the
boundary condition, {the numerical solution of the PDE will
  require discretization with $\cO(M^q)$, $q=\frac{d+1}{d}$, degrees
  of freedom} and, assuming direct solution of the resulting
algebraic problems, the cost of evaluating $z$ on $\Omega$ will be
$\cO(M^{3q})$.  {Thus,} for simplicity, we will restrict our
attention to problems in which the cost of approximating $z$ on $M$
points/elements discretizing the domain $\Omega$ will be $C_1(M^{3q}+M)$,
$q\in\mN$ ($q=0$ represents the case when the dependence of
$e$ on $u$ does not involve a PDE).

A similar argument applies to the evaluation of the second derivative
$e''(u;v,w) = \int_{\Omega} z(v) w \, dx$, except that now $z = z(v)$.
As the {operator defining the} adjoint PDE is the same for
{both} $z$ and $z(v)$, {to determine} $z(v)$ we {only
  need} to perform a {back-}substitution at a computation cost
$C_2 M^{2q}$, {as explained below}.

We note that the cost of evaluating the gradient $g$ corresponding to
a certain $\lambda$ (or equivalently $\eta$) and its derivatives
$g'_i$, $g''_{ij}$, see \eqref{eq:g(eta)}, \eqref{eq:g'i},
\eqref{eq:g''}, will primarily depend on $M$. In general, solution of
each problem of this type requires $\cO(M^3)$ operations.  However,
when several such problems need to be solved with the same
differential operator, then it is more efficient to perform an LU-type
matrix factorization, at the cost $C_3 M^3$, followed by solution of
individual problems via back-substitution, each at the cost $C_4 M^2$.

With these estimates in place and assuming $K \ll M$ and $N_g, N_n
\ll M$, we are now in the position to characterize the complexity of
Algorithms \ref{alg:etagrad} and \ref{alg:etanewton}.  The cost of a
single iteration of the gradient-minimization approach in Algorithm
\ref{alg:etagrad} will be dominated by:
\begin{itemize}
\item[g.1)] one evaluation of $z$ at the computational cost $C_1
  (M^{3q}+M)$,

\item[g.2)] one evaluation of $z(g)$ at the computational cost $C_2
  M^{2q}$,

\item[g.3)]  the following computations {repeated $N_g$ times:}
\begin{itemize}
\item[i.1)] $N+K$ elliptic solves (with factorization) for $g$, $g'_i$
  and $g(\eta-\tau F(\eta))$ at the cost $C_3 M^3 + C_4 (N+K) M^2$,

\item[i.2)] $N+K$ evaluations of $e'(u;v)$ ($v=g$, $v=g'_i$,
  $v=g(\eta-\tau F(\eta))$), and $N+K$ evaluations of $e''(u;g,v)$
  ($v=g$, $v=g'_i$) at the cost at the cost $C_5 (N+K) M$.
\end{itemize}
\end{itemize}
Thus, finding $\hlambda$ and $g(\hlambda)$ with Algorithm
\ref{alg:etagrad} will require
\begin{equation}
{\cal C}_g
=
\cO(1)
\left(
M^{3q} 
+
(M^3 + (N+K)M^2)N_g
\right)
\approx
\cO(1)
\left(
M^{3q} 
+
M^3N_g
\right) 
\
\mbox{\it flops}.
\label{eq:cost(gradient)}
\end{equation}

The cost of a single iteration of Newton's approach in Algorithm
\ref{alg:etanewton} will be dominated by:
\begin{itemize}
\item[n.1)] 
one evaluation of $z$ at the computational cost  $C_1 M^{3q}$,

\item[n.2)] 
$N$ evaluation of $z(g'_i)$ at the computational cost 
$C_2 N M^{2q}$,
\item[n.3)] the following computations {repeated $N_n$ times:}
\begin{itemize}
\item[i.1)] $\frac{1}{2}N^2$ elliptic solves (with factorization) for
  $g''_{ij}$ (noting that $g''_{ij} = g''_{ji}$) at the total cost
  proportional to $C_3 M^3 + C_4 N^2 M^2$,

\item[i.2)] $\frac{1}{2}N^2$ evaluations of $e'(u;v)$ (with $v=g''_{ij}$,
  $1\leq i\leq j\leq N$), and
  $\frac{1}{2}N^2$ evaluations of $e''(u;v,w)$ (with $(v,w)=(g,g''_{ij})$,
  $(v,w)=(g'_i,g'_j)$, $i,j=1,\dots,N$),   
  at the cost $C_6 N^2 M$,

\item[i.3)] one evaluation of $[\partial_j F_i(\eta)]^{-1}\cdot
  [F_i(\eta)]$ at the cost $C_7 N^3$.
\end{itemize}
\end{itemize}
Thus, the cost for computing $\hlambda$ and $g(\hlambda)$ with
Algorithm \ref{alg:etanewton} would require 
\begin{eqnarray}
{\cal C}_{n}
&=&
\cO(1)
\left(
M^{3q}+ N M^{2q}
+
(M^3 + N^2 M^2 +N^3)N_n
\right)
\nonumber\\
&\approx&
\cO(1)
\left(
M^{3q}+ N M^{2q} +
(M^3 + N^2 M^2)N_n
\right)
    \quad
    \mbox{\it flops}. 
\label{eq:cost(newton)}
\end{eqnarray}
Note that the cost of an iteration of a simple gradient algorithm is
\begin{equation}
{\cal C}_{sg}
=
\cO(1)
\left(
M^{3q}+ M^3
\right)
\quad
\mbox{\it flops},
\label{eq:cost(g)}
\end{equation}
Then we {obtain}
\begin{alignat}{2}
 \lim_{N/M\to0}
 \frac{{\cal C}_g}{{\cal C}_{sg}}
 & =
 {\cal O}(1)
  \left(1+M^{3(1-q)}N_g\right),
 & \quad
 \lim_{N/M\to1}
 \frac{{\cal C}_g}{{\cal C}_{sg}}
 & =
 {\cal O}(1)
  \left(1+M^{3(1-q)}N_g\right),
\label{eq:Cg/Csg} 
\\
 \lim_{N/M\to0}
 \frac{{\cal C}_{n}}{{\cal C}_g}
 & =
 {\cal O}(1)
 \frac{1+M^{3(1-q)}N_n}{1+M^{3(1-q)}N_g},
& \quad
 \lim_{N/M\to1}
 \frac{{\cal C}_{n}}{{\cal C}_g}
 & =
 {\cal O}(1)
 \frac{1+MM^{3(1-q)}N_n}{1+M^{3(1-q)}N_g}.
\label{eq:Cn/Cg}
\end{alignat}
Equations \eqref{eq:Cg/Csg} show that the ratio of the cost of our
method using Algorithm \ref{alg:etagrad} {and} the cost of the simple
gradient method is of the same order ${\cal
  O}(1)\left(1+M^{3(1-q)}N_g\right)$, regardless {of} $N$.
Furthermore, the methods tend to {have a comparable cost when $q\geq1$
  and $M$ is} large.  In view of \eqref{eq:Cn/Cg}, it follows that the
same conclusion {also} holds when comparing our method using Algorithm
\ref{alg:etagrad} and Algorithm \ref{alg:etanewton} for $N \ll M$.
{However,} when $N\approx M$, equations \eqref{eq:Cn/Cg} {indicate}
that the cost of our method with Algorithm \ref{alg:etanewton} becomes
substantially higher {(by a factor of $M$) as} compared to the cost
when Algorithm \ref{alg:etagrad} {is used}.  These comments suggest
that it {may be more cost efficient} to use Algorithm
\ref{alg:etagrad} with {large $N$} (under the assumption $K \ll M$),
or Algorithm \ref{alg:etanewton} with $N \ll M$.  {In either case,}
the cost will depend also on $N_g$ {and} $N_n$, i.e., on how fast
Algorithms \ref{alg:etagrad} {and \ref{alg:etanewton} can converge to
  $\heta$}.  In conclusion, the relative efficiency of {original
  Algorithm \ref{alg:gradient} versus its versions using Algorithms
  \ref{alg:etagrad} or \ref{alg:etanewton} to find the optimal
  gradients} will depend on the extend to which the increased
per-iteration cost {in the latter cases} can be offset by the reduced
number of iterations.  This trade-off is illustrated based on a simple
model in the next section.

\section{A model problem and computational results}
\label{sec:results}

In order to illustrate the approach developed in this study, in the
present section we consider the following model problem defined on the
domain $\Omega = (-1,1)$
\begin{equation}
 e(\hu)
 =
 \inf\left\{e(u)
 :=
 \int_\Omega 
\left( 1+ a \, u^2 + a \, \left( \frac{du}{dx} \right)^2 \right)^{1/2} \, dx,
\quad u \in H^1_0(\Omega) \right\},
\label{eq:E}
\end{equation}
where $a = a(x) = 1 - x^2 / 2$. Clearly, the solution is $\hu=0$ and
$e(\hu) = 2$. Energy \eqref{eq:E} gives rise to the following
expressions for its first and second derivative
\begin{align*}
 e'(u;v) 
& = \int_{-1}^1 \left\{ \frac{ a u  }
{\left[ 1+ a \, u^2 + a \, \left( \frac{du}{dx} \right)^2 \right]^{1/2}} - 
\frac{d}{dx} \left(\frac{ a \frac{du}{dx}  }
{\left[ 1+ a \, u^2 + a \, \left( \frac{du}{dx} \right)^2 \right]^{1/2}}\right)\right\} v \, dx, \\
 e''(u;v,w) & = \int_{-1}^1 \left\{ \frac{ a v w  + a \frac{dv}{dx} \frac{dw}{dx} }
{\left[ 1+ a \, u^2 + a \, \left( \frac{du}{dx} \right)^2 \right]^{1/2}} -
\frac{ \left(a u v  + a \frac{du}{dx} \frac{dv}{dx}\right)  
\left(a u w  + a \frac{du}{dx} \frac{dw}{dx}\right)}
{\left[ 1+ a \, u^2 + a \, \left( \frac{du}{dx} \right)^2 \right]^{3/2}} \right\} \, dx.
\end{align*}
To solve problem \eqref{eq:E} we will use the initial guess $u_0(x) =
(1-x^2) \cos(6x) e^x$ chosen such that $u_0\in H^1_0(\Omega)$ and it
has a large $H^1$ norm ensuring that $u_0$ is a ``significant
distance'' away from the solution $\hu$.

In order to mimic the setting with a refined discretization of the
domain $\Omega$, i.e., the case when $M \rightarrow \infty$, in our
computations all functions defined on $\Omega$ (i.e, $u$,
$\widetilde{u}$, $\lambda$, $g_0(u)$, $g(\lambda)$, $g'_i(\lambda)$
and $g''_{ij}(\lambda)$) will be approximated using {\tt Chebfun}
\cite{chebfun}. In this approach all the functions involved are
represented in terms of Chebyshev-series expansions truncated
adaptively to ensure that the truncation errors do not exceed a
prescribed bound (typically related to the machine precision). {\tt
  Chebfun} also makes it possible to solve elliptic boundary-value
problems such as \eqref{eq:g(eta)},\eqref{eq:g'i} and \eqref{eq:g''}
with comparable accuracy. By minimizing the errors related to the
discretization in space, this approach allows us to focus on the
effect of the main parameter in the problem, namely, the dimension $N$
of the space $\Lambda$ in which the optimal weights are constructed,
cf.~\eqref{eq:etaN}. In terms of the basis $\{ \ell^i \}_{i=1}^N$ we
take the standard piecewise-linear ``hat'' functions which, unless
stated otherwise, are constructed based on an equispaced grid. With
such data and choice of the discretization parameters, minimization
problem \eqref{eq:E} is already rich enough to reveal the effect of
the parameter $N$ on convergence and the differences between different
approaches.

We now move on to present computational results obtained solving
problem  \eqref{eq:E} using the following approaches:
\begin{itemize}
\item[(a)] steepest-descent method from Algorithm \ref{alg:gradient}
  with Sobolev gradients $g(\lambda_0)$ defined through the inner
  product \eqref{e:(u,v)_lam0} with {\em constant} weight $\lambda_0 =
  10$ (this value of $\lambda_0$ was found by trial and error to
  produce fastest convergence),

\item[(b)] steepest descent method from Algorithm \ref{alg:gradient}
  with optimal Sobolev gradients $g(\hlambda)$ determined using
  Algorithm \ref{alg:etagrad} for different values $N$; at every iteration 
  Algorithm \ref{alg:etagrad} is restarted
  with the same initial guess $\lambda(x) = \lambda_0$,

\item[(c)] Newton's method from Algorithm \ref{alg:newton},

\item[(d)] steepest descent method from Algorithm \ref{alg:gradient}
  with optimal Sobolev gradients $g(\hlambda)$ determined using a
  simplified version of Algorithm \ref{alg:etanewton} for different
  values $N$ (see below for details); at every iteration simplified
  Algorithm \ref{alg:etanewton} is restarted with the same initial
  guess $\lambda(x) = \lambda_0$.
\end{itemize}
Approaches (a), (b) and (d) use the same fixed step size $\kappa =
50$.  Approximations of the exact solution $\widehat{u}$ obtained at
the $n$th iteration will be denoted $u_n$. In order to prevent the
optimal weights $\hlambda(x)$ from becoming too close to zero for
certain $x$, which would complicate the numerical solution of problems
\eqref{eq:g(eta)}, \eqref{eq:g'i} and \eqref{eq:g''}, the line-search
in Algorithm \ref{alg:etagrad} and the length of Newton's step in
Algorithm \ref{alg:etanewton} are restricted such that $\min_{x \in
  [-1,1]} \hlambda(x) > \epsilon_{\tau} \lambda_0$, where we used
$\epsilon_{\tau} = 10^{-2}$. In addition, since this will make it
possible to objectively compare cases with different values of $N$,
here we modify the termination condition in Algorithm
\ref{alg:etagrad}, cf.~line 13, by replacing it with one given in
terms of a minimum relative decrease of $j(\eta)$, i.e., $|j(\heta) -
j(\eta)| / j(\eta) \le \epsilon_{\lambda}$, where $\epsilon_{\lambda}$
is a prescribed tolerance.

We now examine the effect of different parameters on the results
obtained with each of the approaches (a)--(d) defined above.

\subsection{Analysis of the effect of the tolerance $\epsilon_\lambda$}
\begin{figure}[!h]
\centering
\mbox{
\subfloat[]{\includegraphics[width=0.475\textwidth]{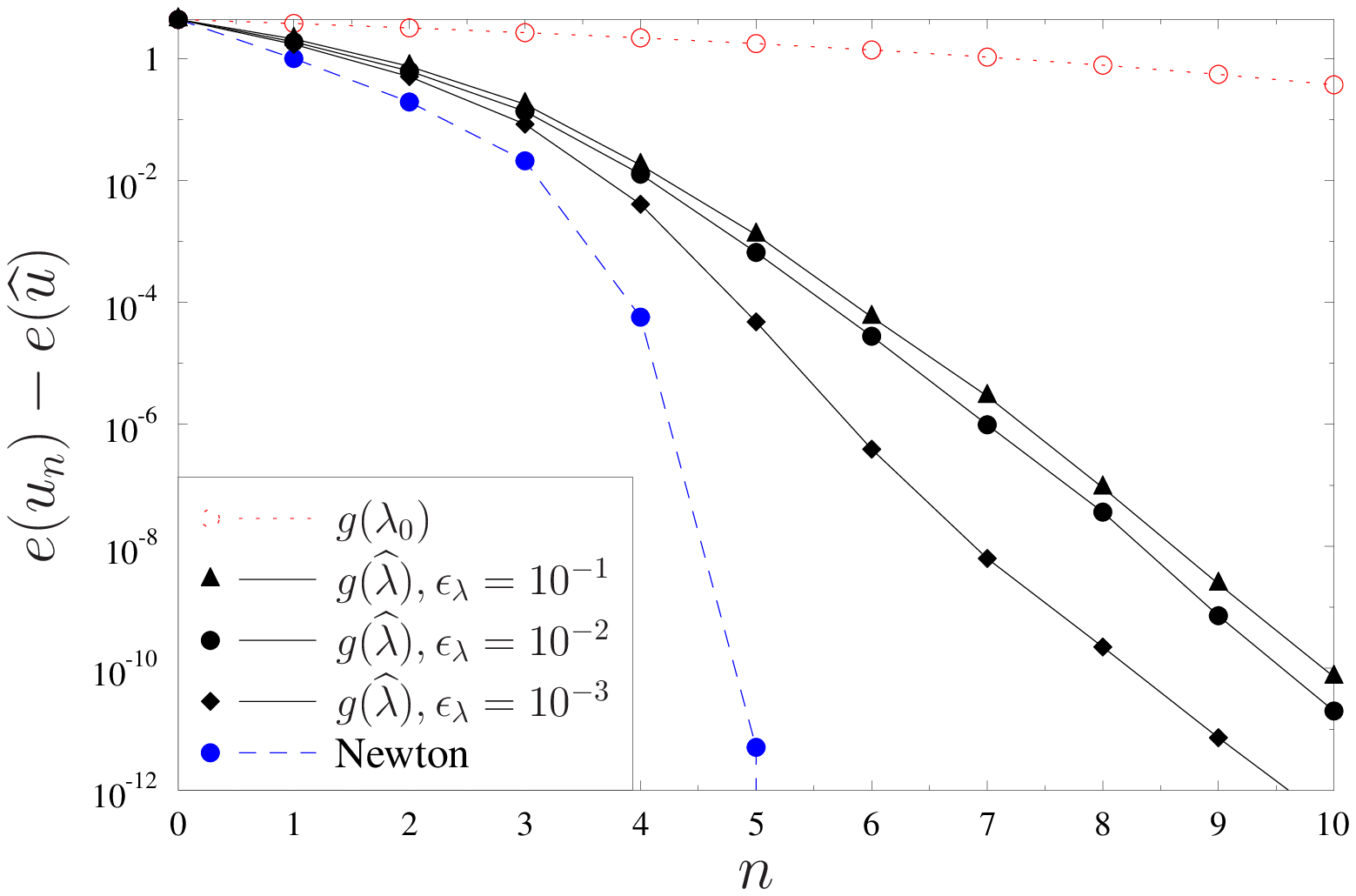}}
\quad
\subfloat[]{\includegraphics[width=0.475\textwidth]{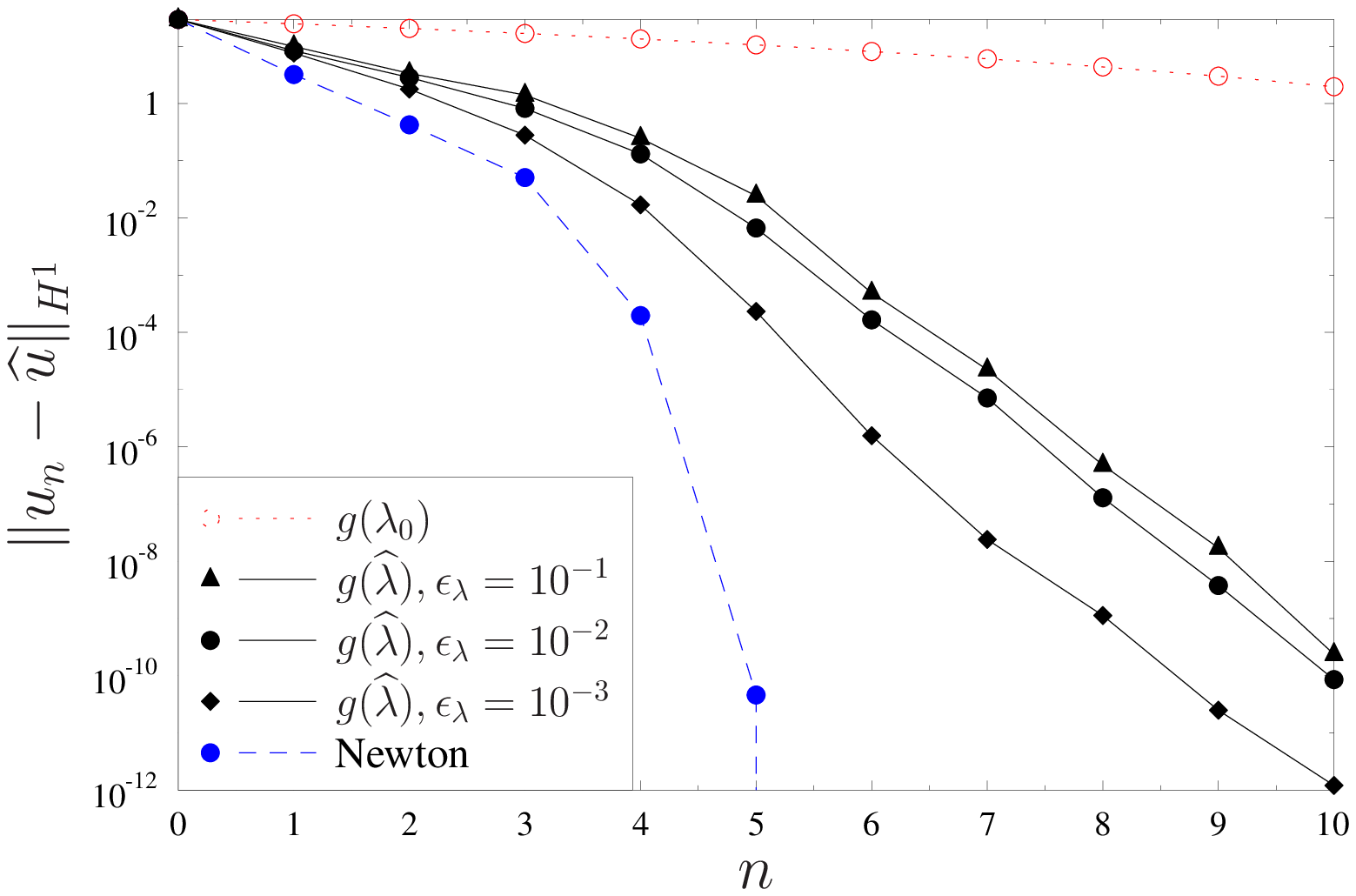}}
}
\caption{(a) (Shifted) energy $e(u_n)-e(\hu)$ and (b) the $H^1$
  approximation error $\|u_n - \hu \|_{H^1}$ as functions of the
  iteration count $n$ where the optimal gradients $g(\hlambda)$ are
  obtained {using} Algorithm \ref{alg:etagrad} with $N=50$ and
  different tolerances $\epsilon_{\lambda} = 10^{-1}, 10^{-2},
  10^{-3}$. For comparison, the results obtained using standard
  Sobolev gradients $g(\lambda_0)$ with a constant weight $\lambda_0 =
  10$ and with Newton's method, cf.~Algorithm \ref{alg:newton}, are
  also presented.}
\label{fig:eps}
\end{figure}

The decrease of the (shifted) energy $e(u_n)-e(\hat{u})$ and of the
$H^1$ approximation error $\|u_n - \hu\|_{H^1}$ are shown for
approaches (a), (b) and (c) in Figures \ref{fig:eps}a and
\ref{fig:eps}b, respectively, where in case (b) we used a single value
$N=50$ and three different tolerances $\epsilon_{\lambda} = 10^{-1},
10^{-2}, 10^{-3}$.  In Figure \ref{fig:eps}a we see that minimization
with optimal gradients $g(\hlambda)$ produces a significantly faster
decrease of energy $e(u_n)$ than optimization with ``standard''
Sobolev gradients $g(\lambda_0)$ and analogous trends are also evident
in the decrease of the approximation error $\|u_n - \hu\|_{H^1}$,
cf.~Figure \ref{fig:eps}b. We add that in order to solve the
minimization problem to the same level of accuracy the method based on
the ``standard'' Sobolev gradients $g(\lambda_0)$ requires as many as
42 iterations (for clarity, these later stages are not shown in the
figures). 

In addition, in Figures \ref{fig:eps}a and \ref{fig:eps}b
we also observe that convergence of the proposed method systematically
accelerates as the tolerance $\epsilon_{\lambda}$ is refined, i.e., as
the optimal weights $\hlambda$ are approximated more accurately.
However, we remark that reducing $\epsilon_{\lambda}$ below $10^{-3}$
did not produce further improvement of convergence.  Hence, hereafter
we will set $\epsilon_{\lambda} = 10^{-3}$.

\subsection{Analysis of the effect of the dimension $N$ of the
approximation space $\Lambda$}
\begin{figure}[!h]
\centering
\mbox{
\subfloat[]{\includegraphics[width=0.475\textwidth]{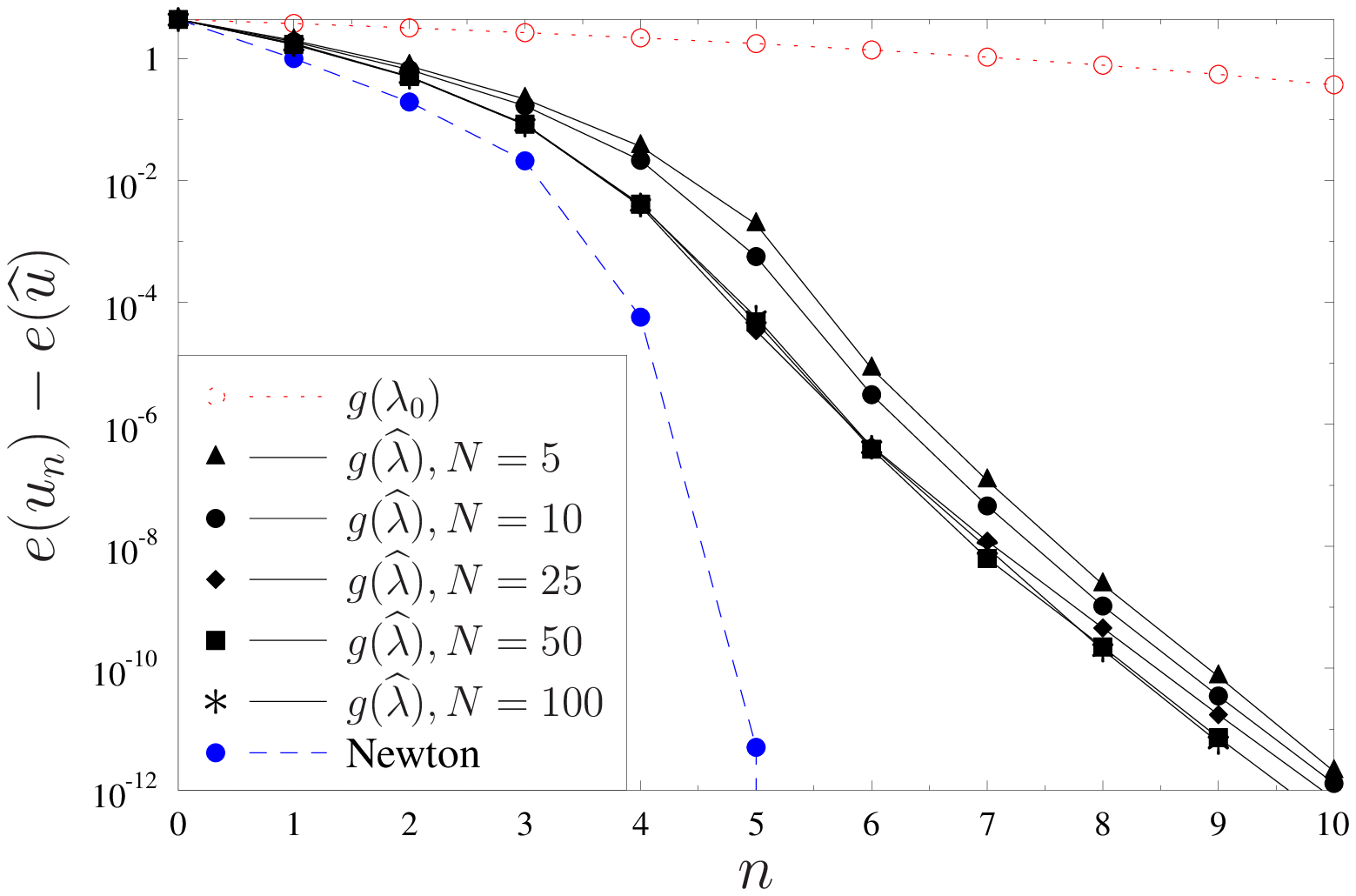}}
\quad
\subfloat[]{\includegraphics[width=0.475\textwidth]{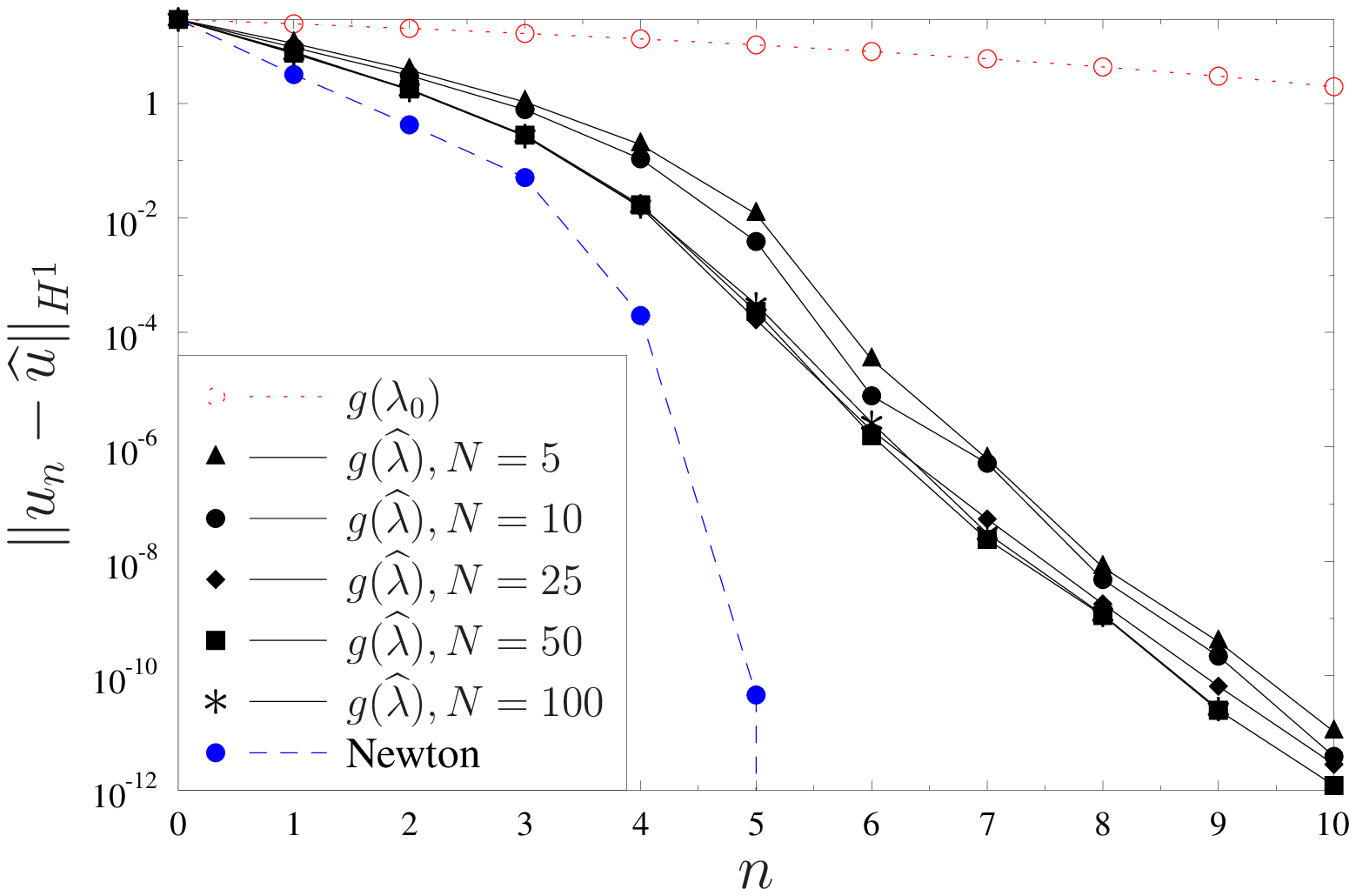}}
}
\caption{(a) (Shifted) energy $e(u_n)-e(\hu)$ and (b) the $H^1$
  approximation error $\|u_n - \hu \|_{H^1}$ as functions of the
  iteration count $n$ for different dimensions $N$ where the optimal
  gradients $g(\hlambda)$ are obtained with Algorithm
  \ref{alg:etagrad}. For comparison, the results obtained using
  standard Sobolev gradients $g(\lambda_0)$ with a constant weight
  $\lambda_0 = 10$ and with Newton's method, cf.~Algorithm
  \ref{alg:newton}, are also presented.}
\label{fig:min}
\centering
\mbox{
\subfloat[]{\includegraphics[width=0.475\textwidth]{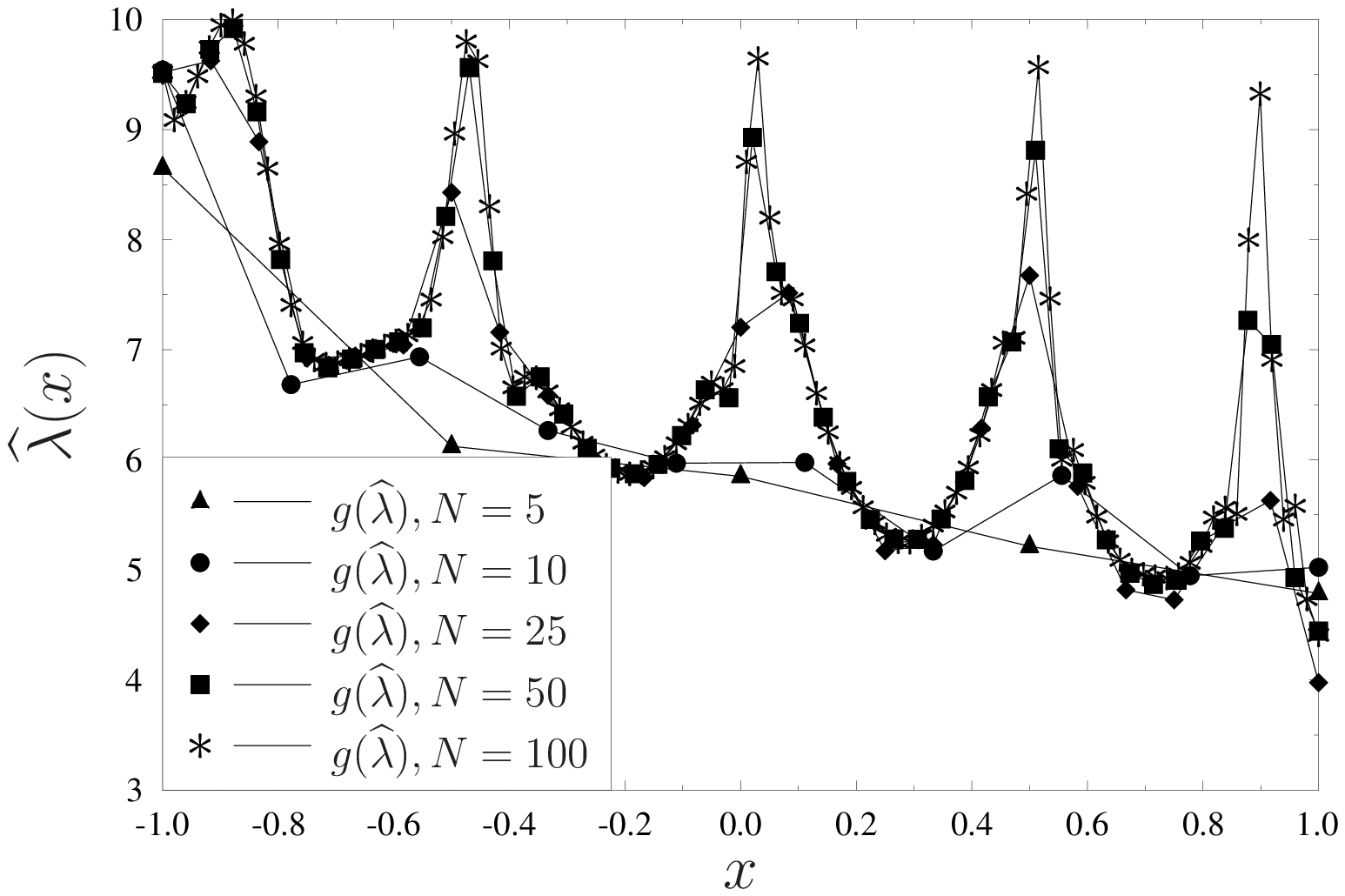}}
\quad
\subfloat[]{\includegraphics[width=0.475\textwidth]{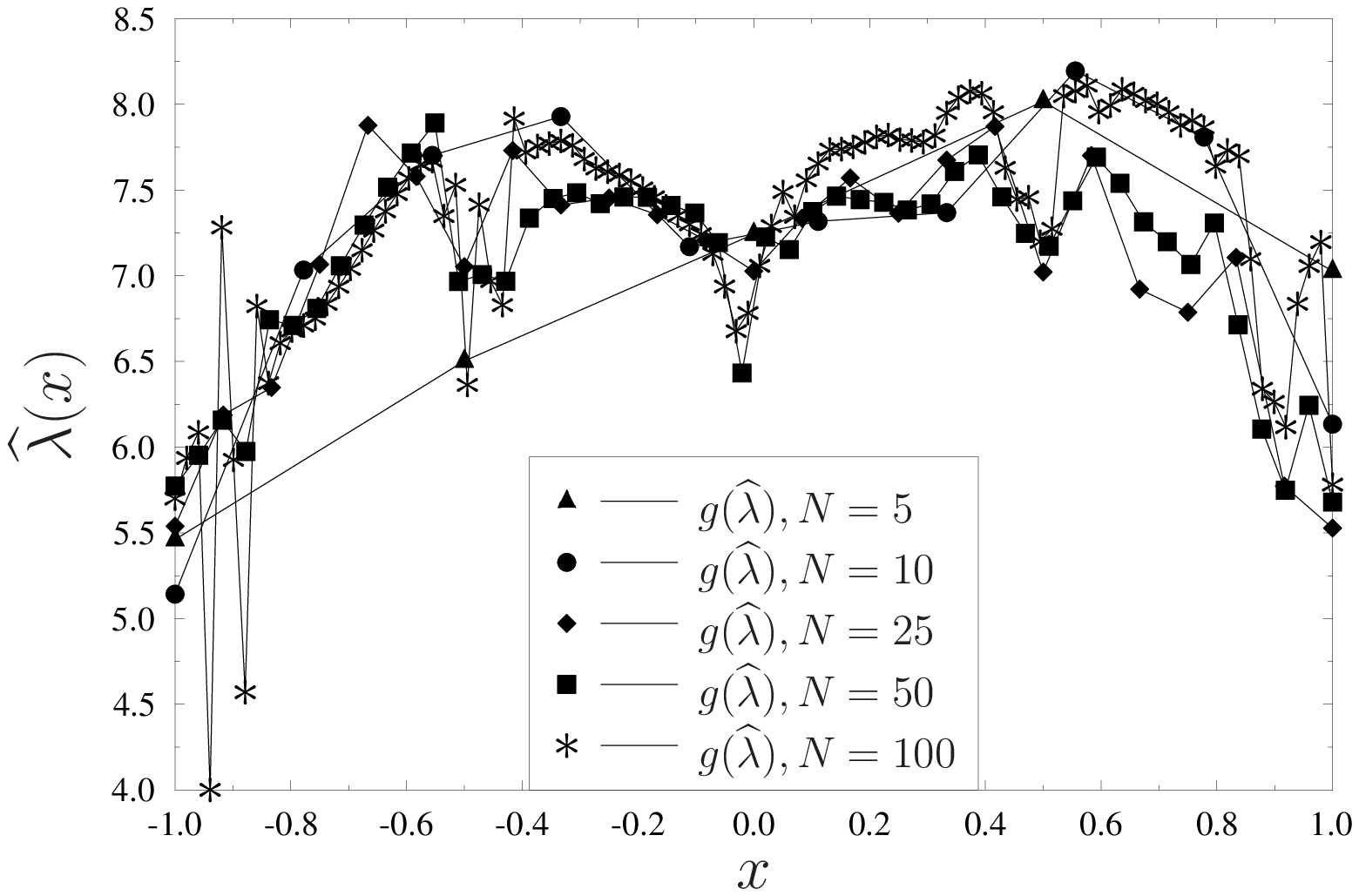}}
}
\caption{Optimal weights $\hlambda(x)$ as functions of $x$ obtained
  for different $N$: {(a) an early ($n=2$) iteration and (b) a
    late ($n=8$)} iteration of Algorithm \ref{alg:etagrad}, cf.~Figure
  \ref{fig:min}.}
\label{fig:lam}
\end{figure}

{The results concerning} the effect of $N$ on the performance of
approach (b) are compared with the data for approaches (a) and (c) in
Figures \ref{fig:min}a and \ref{fig:min}b for the (shifted) energy
$e(u_n)-e(\widehat{u})$ and the $H^1$ approximation error $\|u_n -
\hu\|_{H^1}$, respectively. We observe that, when optimal gradients
$g(\hlambda)$ are used, both $e(u_n)-e(\hu)$ and $\|u_n - \hu\|_{H^1}$
initially reveal a quadratic convergence, similar to the behavior of
these quantities in Newton's method, followed at later iterations by a
linear convergence, typical of the standard gradient method. In the
light of Theorem \ref{th:P+newton->1}, cf.~estimate
\eqref{e:|P(tu-u)|->1}, this observation can be explained by the fact
that at early iterations dominant components of the error $(u_n-\hu)$
are contained in the subspaces $T_{u_n}$ where the optimal gradients
$g(\hlambda)$ are consistent with Newton's steps $h$, cf.~Remark
\ref{r:name}. Then, once these error components are eliminated, at
later iterations the error $(u_n - \hu)$ is dominated by components in
directions orthogonal to $T_{u_n}$ where the optimal gradients
$g(\hlambda)$ do not well reproduce the Newton steps $h$.  In Figures
\ref{fig:min}a and \ref{fig:min}b we also see that the convergence
improves as the dimension $N$ is increased until it saturates for $N$
large enough (here $N \gtrapprox 25)$.  This could be explained by the
conjecture that increasing $N$ above a certain limit (approximately
$25$ in this case) does not increase the ``effective'' dimension of
$T_{u_n}$ in $H$ anymore (such a possibility is allowed by the error
analysis presented in Section \ref{s:analysis->1}).

In this context it
is also interesting to investigate the evolution of the spatial
structure of the optimal weights $\hlambda(x)$ and these results are
shown for different values of $N$ at an early ($n=2$) and a later
($n=8$) iteration in Figures \ref{fig:lam}a and \ref{fig:lam}b,
respectively.  In the first case ($n=2$ corresponding to the quadratic
convergence) we see that the optimal weights $\hlambda(x)$ converge to
a well-defined profile as $N$ increases, which features a number of
distinct ``spikes''. On the other hand, at later iterations ($n=8$
corresponding to the linear regime) the convergence of the optimal
weights $\hlambda(x)$ with $N$ is less evident and the resulting
profiles tend to be more uniform.

\begin{figure}[!h]
\centering
\mbox{
\subfloat[]{\includegraphics[width=0.475\textwidth]{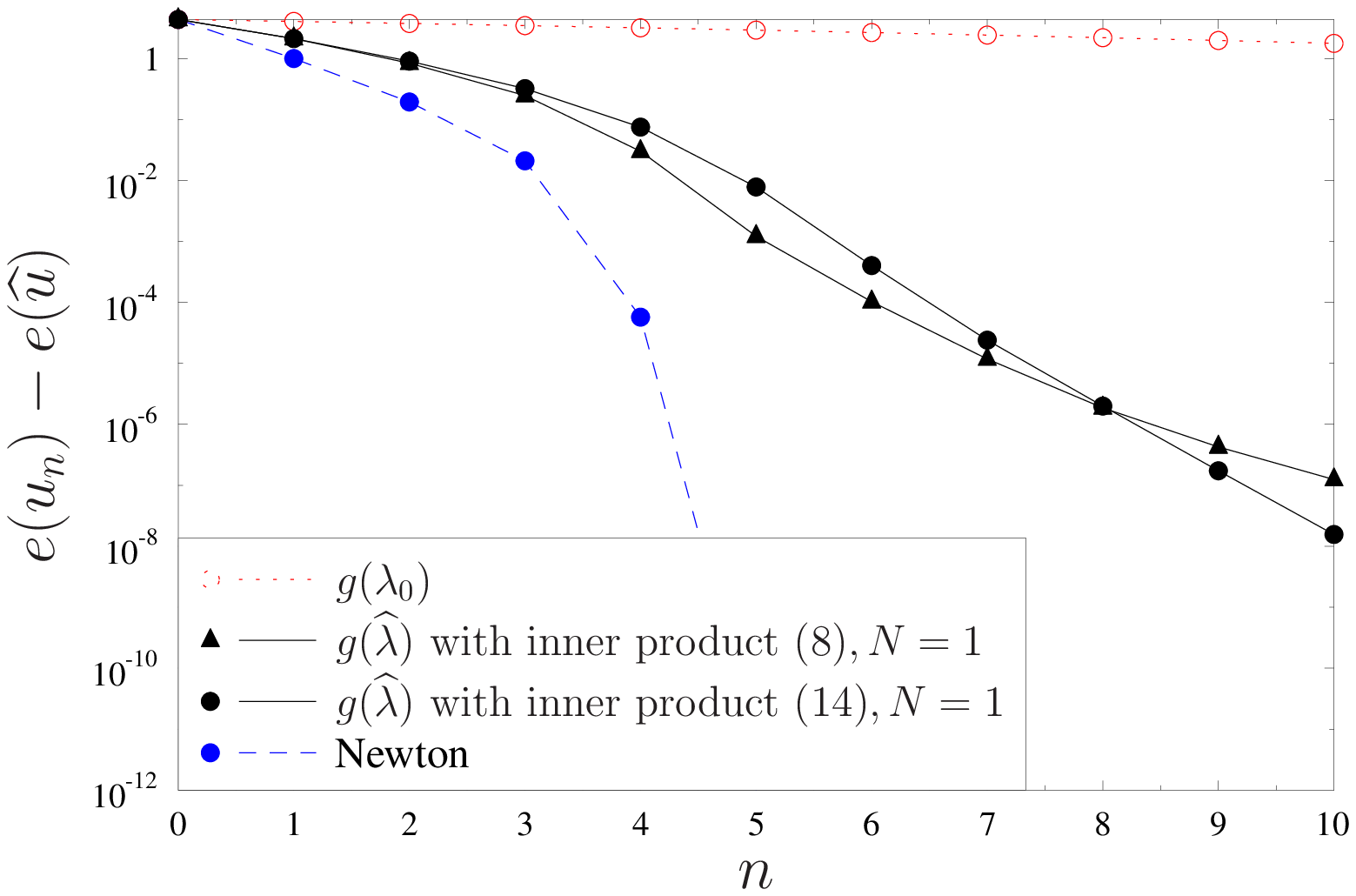}}
\quad
\subfloat[]{\includegraphics[width=0.475\textwidth]{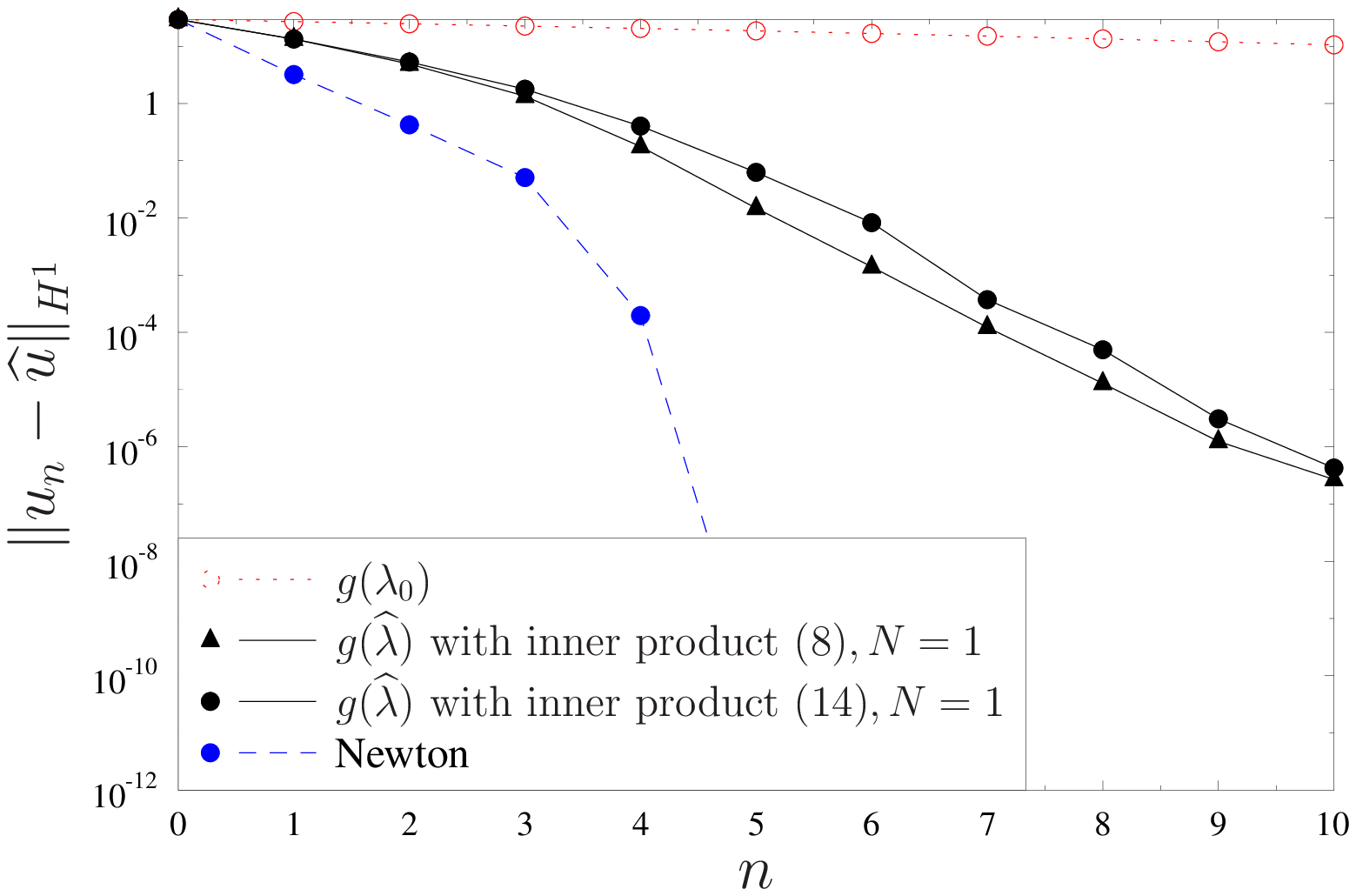}}
}
\caption{(a) (Shifted) energy $e(u_n)-e(\hu)$ and (b) the $H^1$
  approximation error $\|u_n - \hu \|_{H^1}$ as functions of the
  iteration count $n$ for the case when $N=1$ and the optimal
  gradients $g(\hlambda)$ are obtained using the inner product
  definitions \eqref{e:(u,v)_lambda} combined with Algorithm
  \ref{alg:etagrad} and \eqref{eq:ipH10} combined with {the}
  explicit relation \eqref{e:hlambda,N=1}. For comparison, the results
  obtained using standard Sobolev gradients $g(\lambda_0)$ with a
  constant weight $\lambda_0 = 10$ and with Newton's method,
  cf.~Algorithm \ref{alg:newton}, are also presented. The step size
  used in these calculations is $\kappa = 25$.}
\label{fig:N1}
\end{figure}

\noindent
We want to {highlight} the case when $N=1$ and space $H$ is
endowed with the inner product redefined as in \eqref{eq:ipH10}. As
shown in Remark \ref{r:hlambda,H10}, in such circumstances the optimal
$\lambda$ can be found analytically, cf.~relation
\eqref{e:hlambda,N=1}, at essentially no cost and the iterations
produced by Algorithm \ref{alg:gradient} do not depend on the step
size $\kappa$. The results obtained with this approach and using the
optimal gradients $g(\hlambda)$ defined in terms of the inner product
\eqref{e:(u,v)_lambda} are compared in Figures \ref{fig:N1}a and
\ref{fig:N1}b for the (shifted) energy $e(u_n)-e(\widehat{u})$ and the
$H^1$ approximation error $\|u_n - \hu\|_{H^1}$, respectively. As is
evident from these figures, the performance of the approaches
corresponding to the two definitions of the inner product,
\eqref{e:(u,v)_lambda} and \eqref{eq:ipH10}, is comparable and
{in both cases} much better than when a fixed weight $\lambda_0$
is used. We stress that in the case corresponding to the inner product
\eqref{eq:ipH10} determination of the optimal $\lambda$ does not
require an iterative solution.


\subsection{Analysis of the robustness of approach (b) with respect to
  {variations of the basis functions defining $\eta$}}
\begin{figure}[H]
\centering
\mbox{
\subfloat[]{\includegraphics[width=0.475\textwidth]{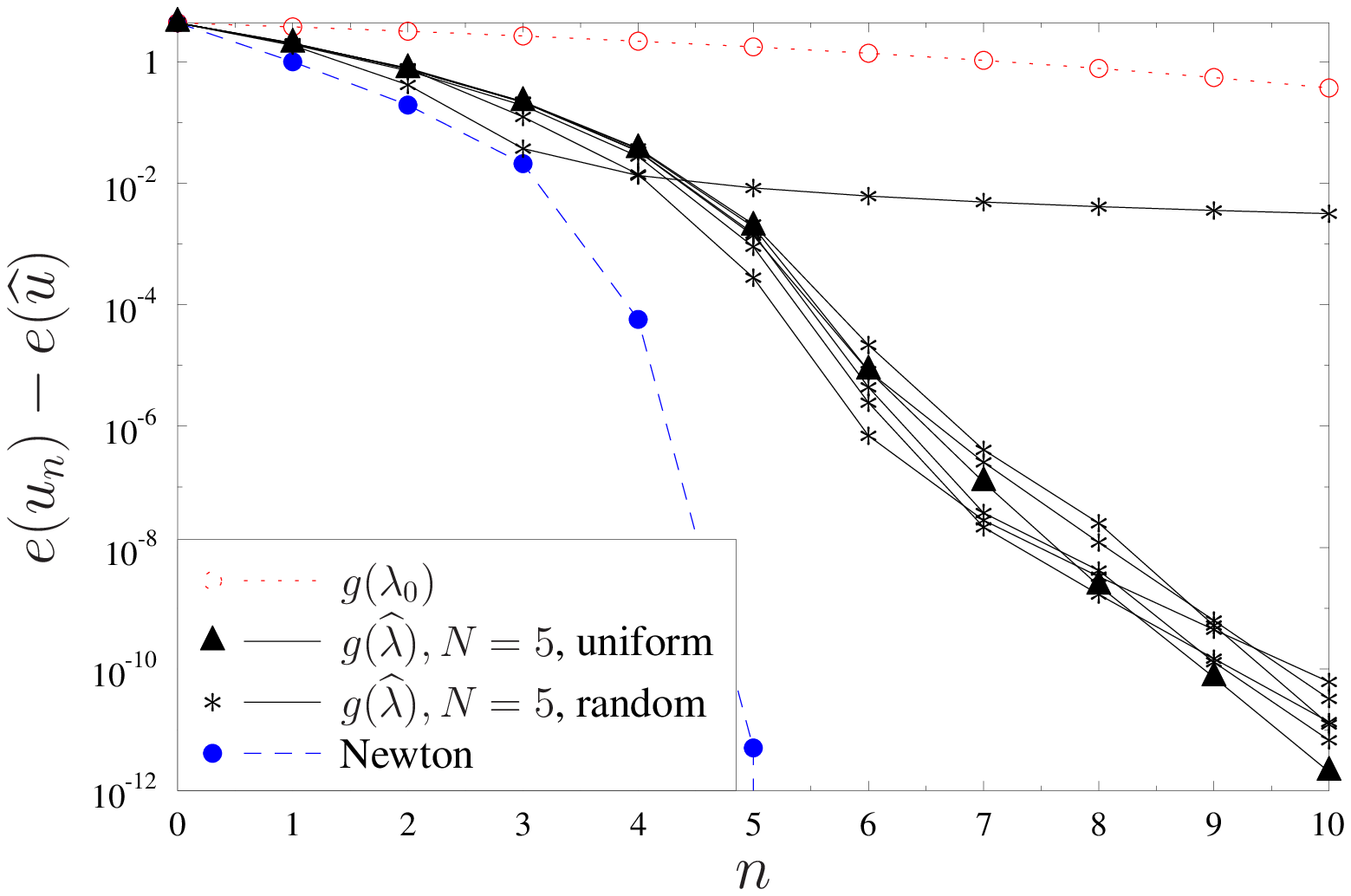}}
\quad
\subfloat[]{\includegraphics[width=0.475\textwidth]{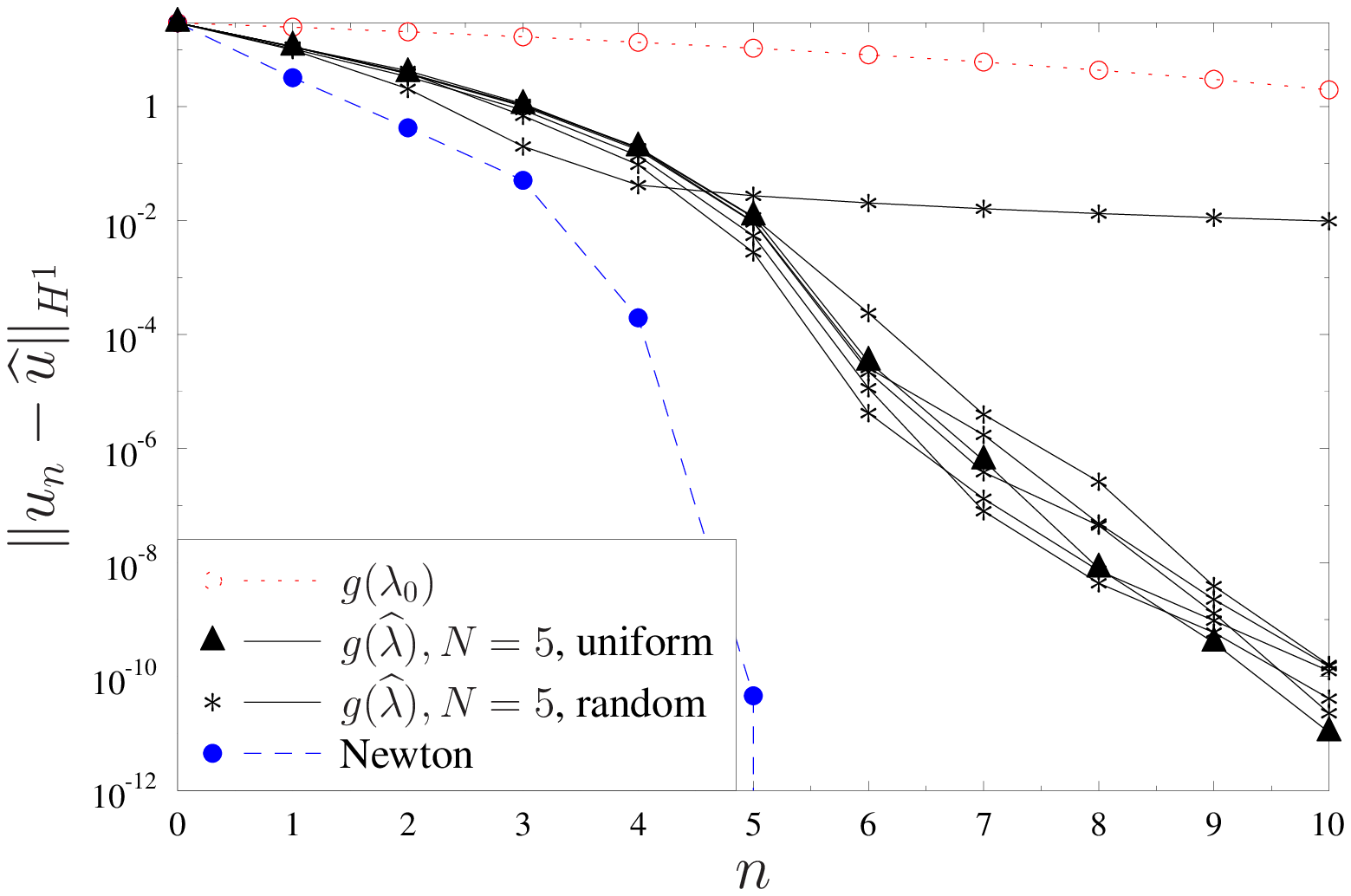}}
}
\caption{(a) (Shifted) energy $e(u_n)-e(\hu)$ and (b) the $H^1$
  approximation error $\|u_n - \hu \|_{H^1}$ as functions of the
  iteration count $n$ for the case when $N=5$ and the optimal
  gradients $g(\hlambda)$ are obtained with Algorithm
  \ref{alg:etagrad} using uniform and random distributions of grid
  points defining the basis functions $\{ \ell^i \}_{i=1}^N$. For
  comparison, the results obtained using standard Sobolev gradients
  $g(\lambda_0)$ with a constant weight $\lambda_0 = 10$ and with
  Newton's method, cf.~Algorithm \ref{alg:newton}, are also
  presented.}
\label{fig:rand}
\end{figure}

This {analysis is performed} by constructing basis functions $\{
\ell^i \}_{i=1}^N$ based on grid points distributed randomly with an
uniform probability distribution over the interval $(-1,1)$, except
for the leftmost and the rightmost grid points which are always at $x
= \pm 1$. The results obtained in several realizations with $N=5$ are
compared to the reference case of basis functions constructed based on
equispaced grid points as well as with the results obtained with
approaches (a) and (c) in Figures \ref{fig:rand}a and \ref{fig:rand}b
for the (shifted) energy $e(u_n)-e(\widehat{u})$ and the $H^1$
approximation error $\|u_n - \hu\|_{H^1}$, respectively. One can see
in these figures that, expect for one {realization} corresponding
to a very special distribution of the grid points, the convergence is
little affected by the choice of the basis $\{ \ell^i \}_{i=1}^N$.

\subsection{{Analysis of the performance of a simplified version of Algorithm \ref{alg:etanewton}}}
\begin{figure}[!h]
\centering
\mbox{
\subfloat[]{\includegraphics[width=0.475\textwidth]{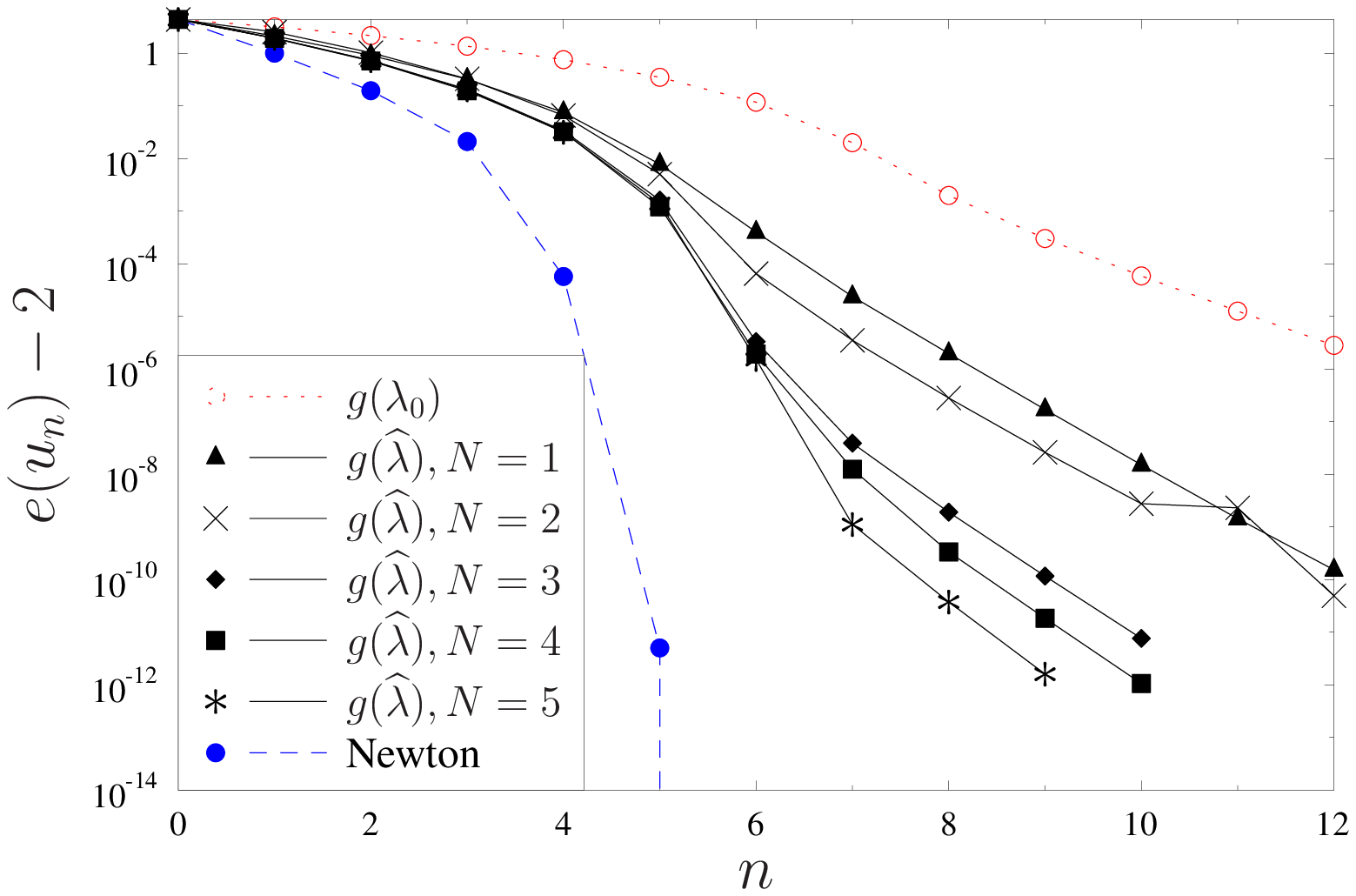}}
\quad
\subfloat[]{\includegraphics[width=0.475\textwidth]{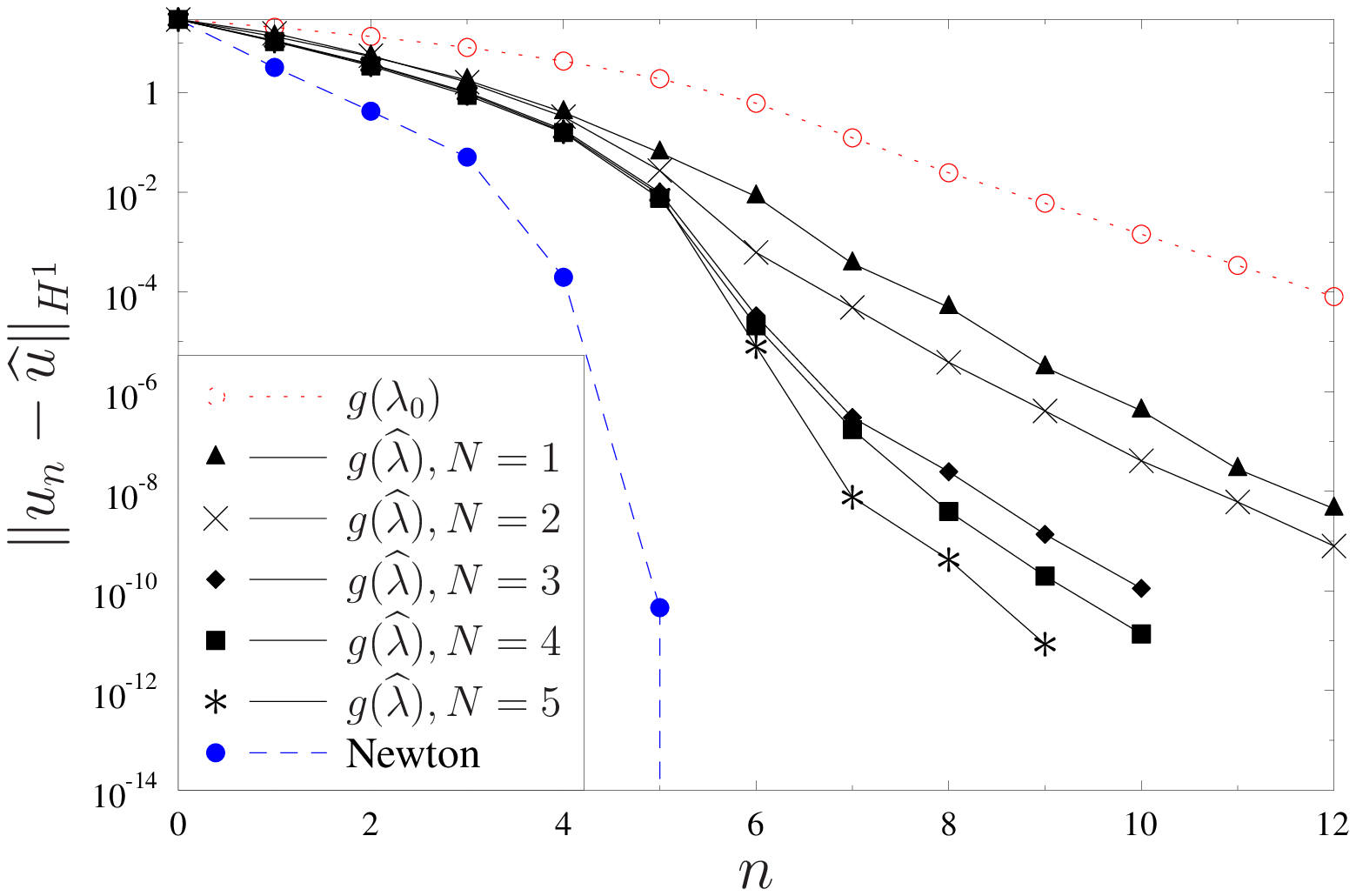}}
}
\caption{(a) (Shifted) energy $e(u_n)-e(\hu)$ and (b) the $H^1$
  approximation error $\|u_n - \hu \|_{H^1}$ as functions of the
  iteration count $n$ for different dimensions $N$ where the optimal
  gradients $g(\hlambda)$ are obtained with the simplified version of
  Algorithm \ref{alg:etanewton} described above.  For comparison, the
  results obtained using standard Sobolev gradients $g(\lambda_0)$
  with a constant weight $\lambda_0 = 10$ and with Newton's method,
  cf.~Algorithm \ref{alg:newton}, are also presented.}
\label{fig:lin}
\end{figure}

Finally, we consider approach (d) where the optimal weights
$\hlambda(x)$ and the corresponding optimal gradients $g(\hlambda)$
are determined with Algorithm \ref{alg:etanewton} simplified as
follows. The complexity analysis presented in Section
\ref{sec:complexity} shows that Algorithm \ref{alg:etanewton} may be
quite costly from the computational point of view when $N \gg 1$. To
alleviate this difficulty, we consider its simplified version where
only one iteration ($N_n = 1$) is performed on system
\eqref{eq:alpha->F=0} in which the ``test'' functions $g'_i$,
$i=1,\dots,N$, are assumed not to depend on $\lambda$ (or $\eta$). In
other words, since instead of $g'_i(\eta)$, $i=1,\dots,N$, the
functions $g'_i(\eta_0)$ are used to obtain expressions for
$[F(\eta)]_i$, $i=1,\dots,N$, in \eqref{eq:alpha->F=0}, the second
derivatives $g''_{ij}$ are eliminated from the Hessian
$[DF(\eta)]_{ij}$, $i,j=1,\dots,N$ in \eqref{eq:DF}, which very
significantly reduces the computational cost.  The results obtained
with this simplified approach are shown in Figures \ref{fig:lin}a and
\ref{fig:lin}b, respectively, for the decrease of the (shifted) energy
$e(u_n)-e(\hu)$ and for the decrease of the $H^1$ approximation error
$\|u_n - \hu\|_{H^1}$. In these figures we observe general trends
qualitatively similar to those evident in Figures \ref{fig:min}a and
\ref{fig:min}b, except that the convergence is slower and the
transition from the quadratic to linear convergence tends to occur at
earlier iterations.

\section{Conclusions}
\label{sec:final}
We have developed a gradient-type numerical approach for
{unconstrained} optimization problems in infinite-dimensional
Hilbert spaces.
Our method consists in
finding an optimal inner product among a family of equivalent inner
products parameterized by a space-dependent weight $\lambda$ function. 
The optimal weight $\hlambda$ solves a nonlinear optimization problem in a 
finite dimensional subspace.
Rigorous analysis demonstrates that, in addition to the linear
convergence characterizing the standard gradient method, the proposed
approach also attains quadratic convergence in the sense that the
projection error in a finite-dimensional subspace generated in the
process decreases quadratically. Or, equivalently, in this finite
dimensional subspace, the optimal gradients and Newton's steps are
equivalent.  The dimension of these subspaces is determined by the
number $N$ of discrete degrees of freedom parameterizing the inner
products through the weight $\lambda$.

This analysis is confirmed by numerical experiments,
performed based on a simple optimization problem in a setting
mimicking high spatial resolution. More specifically, at early
iterations both the minimized energy and the error with respect to the
exact solution exhibit quadratic convergence followed by the linear
convergence at later iterations. The behavior of the proposed method
also reveals expected trends with variations of the numerical
parameters, namely, the dimension of the space in which the optimal
weights $\hlambda$ are constructed, properties of the basis in this
space and the accuracy with which the inner optimization problems are
solved. In all cases the convergence of the proposed approach is much
faster than obtained using Sobolev gradients with fixed weights. For
the ease of analysis and comparisons we focused on a gradient-descent
method with a fixed step size $\kappa$, but it can be expected that a
similar behavior will also occur when the step size is determined
adaptively through suitable line minimization.

The complexity analysis performed in Section \ref{sec:complexity}
indicates that the per-iteration cost of the proposed approach
and of the standard Sobolev gradient method have the same order if,
for example, the energy depends on a elliptic PDE.  When the dimension
of weight space is $N = 1$ and the inner product does not have the
$L^2$ term, cf.~\eqref{eq:ipH10}, then the optimal weight is given
explicitly, eliminating the need for its numerical determination. In
this particular case the proposed approach has some similarity to the
Barzilai-Borwein algorithm \cite{barzilai-1} and produces iterates
which do not depend on the step size in the gradient method.  The
computational cost of the proposed approach is also significantly
reduced when Algorithm \ref{alg:etanewton} is used in a simplified
form, as described in Section \ref{sec:results}. {We thus
  conclude that the gradient-descent method from Algorithm
  \ref{alg:gradient} combined with Algorithms \ref{alg:etagrad} and
  \ref{alg:etanewton} used to find optimal gradients} are promising
approaches suitable for large-scale optimization problems and
applications to some such problems will be investigated in the near
future.

Our approach based on optimal gradients differs from the family of
quasi-Newton methods in that instead of approximating the full Hessian
using gradients from earlier iterations, see for example \cite{nw00},
it relies on computing the action of the exact \revt{Hessian and
  gradients,} but only on a few judiciously selected directions, and
then matching them by appropriately choosing the inner product.
Consequently, the resulting algebraic problem is of a much smaller
dimension thereby avoiding complications related to poor conditioning
and computational cost.

Finally, we believe that the analysis and results presented here
explain the acceleration of gradient minimization reported for a range
of different problems in \cite{pbh02,p07b,rssl09,r10,ms10,dk10} when
Sobolev gradients with suitable (constant) weights were used.
Moreover, our work also provides a rational and constructive answer to
the open problem of finding an optimal form of the inner product
defining the Sobolev gradients.

\bibliographystyle{plain}

\end{document}